\documentclass[12pt]{article}
\usepackage{fullpage,graphicx,psfrag,amsmath,amsfonts,amsthm, amssymb, verbatim}
\usepackage{cite}
\usepackage{cases}
\usepackage{epstopdf}
\usepackage{mdwlist}
\usepackage{bigstrut,array,multirow,tabularx}
\usepackage[dvipsnames]{xcolor}
\usepackage{subfigure}
\usepackage{dsfont}
\usepackage{bbm,bm}
\usepackage{url}
\usepackage{algorithm,algorithmic}
\usepackage[final]{changes}

\DeclareMathAlphabet{\mathsf}{OML}{cmbr}{m}{it}

\newtheorem{definition}{Definition}
\newtheorem{theorem}{Theorem}
\newtheorem{lemma}{Lemma}

\newtheorem{proposition}{Proposition}

\newtheorem{remark}{Remark}
\newtheorem{property}{Property}
\newtheorem{assumption}{Assumption}

\DeclareMathOperator{\tr}{\mathrm{tr}}
\DeclareMathOperator{\E}{\mathds{E}}

\DeclareMathOperator*{\argmin}{argmin}

\DeclareMathOperator{\prox}{prox}

\newcommand{\N}{\mathbbmss{N}}
\newcommand{\R}{\mathbbmss{R}}

\newcommand{\B}[1]{\mathbf{#1}}

\newcommand{\EX}[1]{\E\left\{{#1}\right\}}

\newcommand{\IM}[1]{\B{I}_{#1}}

\newcommand{\InP}[2]{\left\langle #1,#2 \right \rangle}

\newcommand{\BB}[1]{\ensuremath{\bm{#1}}}

\newcommand{\xv}{\BB{x}}
\newcommand{\gv}{\BB{g}}
\newcommand{\xTd}{\tilde{\BB{x}}}
\newcommand{\vv}{\BB{v}}
\newcommand{\uv}{\BB{u}}
\newcommand{\dv}{\BB{d}}
\newcommand{\zv}{\BB{z}}
\newcommand{\yv}{\BB{y}}
\newcommand{\tv}{\BB{t}}
\newcommand{\av}{\BB{a}}
\newcommand{\bv}{\BB{b}}

\newcommand{\grad}[1]{\nabla #1}

\newcommand{\He}{\BB{H}}

\newcommand{\Am}{\BB{A}}
\newcommand{\Bs}{\mathcal{B}}

\newcommand{\fF}[2]{f_{#1}\left(#2\right)}
\newcommand{\hF}[1]{h\left(#1\right)}

\newcommand{\FF}[1]{F\left(#1\right)}

\newcommand{\Nvar}{\Delta}
\newcommand{\Prox}[3]{{\prox}_{#1}^{#2}\left(#3\right)}

\newcommand{\Lyp}{V}
\newcommand{\ksub}{\kappa_{\mathrm{sub}}}
\newcommand{\kmax}{\kappa_{\mathrm{max}}}
\newcommand{\kavg}{\kappa_{_{\He}}}
\newcommand{\ktil}{\tilde{\kappa}}
\newcommand{\Lavg}{L_{\mathrm{avg}}}
\newcommand{\Lmax}{L_{\mathrm{max}}}
\newcommand{\reg}{\gamma}
\newcommand{\lvr}{\mathsf{s}}
\newcommand{\sdim}{\mathsf{d}}

\newcommand{\eig}{\lambda}

\newcommand{\norm}[1]{\left|{#1}\right|}
\newcommand{\Vnorm}[1]{\left\|{#1}\right\|}
\newcommand{\Mnorm}[2]{\left\Vert{#1}\right\Vert_{#2}}

\definechangesauthor[color=orange]{VM}

\title{Curvature-Exploiting Acceleration of Elastic Net Computations}

\author{
	{Vien V. Mai and  Mikael Johansson}{\thanks{V. V. Mai and M. Johansson are with the Department of Automatic Control, School of Electrical Engineering and Computer Science, Royal Institute of Technology (KTH), SE-100 44 Stockholm, Sweden. Emails: \tt\small\{maivv, mikaelj\}@kth.se.}%
	}       
}

\begin{document}

\maketitle

\begin{abstract} 

This paper introduces an efficient second-order method for solving the elastic net problem.
Its key innovation is a computationally efficient technique for injecting curvature information in the optimization process which admits a strong theoretical performance guarantee.  
 In particular, we show improved run time over popular first-order methods and quantify the speed-up in terms of statistical measures of the data matrix. The improved time complexity is the result of an extensive exploitation of the problem structure and a careful combination of second-order information, variance reduction techniques, and momentum acceleration. Beside theoretical speed-up, experimental results demonstrate great practical performance benefits of  curvature information, especially for ill-conditioned data sets.
\end{abstract}

\section{Introduction}

Lasso\added{, ridge} and \replaced{e}{E}lastic net regression are\deleted{two} fundamental problems in statistics\replaced{ and}{,} machine learning,\deleted{, and optimization} with countless applications in \replaced{science and engineering~\cite{ZH05}.}{signal processing, data analysis,}\deleted{and many other fields~\cite{ZH05}.}\deleted{among many others.} 
\deleted{Given data matrices $\Am\in\R^{n\times d}$ and $\bv\in \R^n$, elastic} Elastic net regression amounts to solving the following convex optimization problem
\begin{align}\label{eq:elastic}
	\underset{\xv\in\R^d}{\mbox{minimize}}\,
	\bigg\{
	\frac{1}{2n}\Vnorm{\Am\xv-\bv}_2^2+\frac{\reg_2}{2}\Vnorm{\xv}_2^2+\reg_1\Vnorm{\xv}_1
	\bigg\},
\end{align}
for given data matrices $\Am\in\R^{n\times d}$ and $\bv\in \R^n$ and regularization parameters  $\reg_1$ and $\reg_2$. \added{Setting $\gamma_1=0$ results in ridge regression,  $\gamma_2=0$ yields lasso and letting $\gamma_1=\gamma_2=0$ reduces the problem to the classical least-squares. Lasso promotes sparsity of the optimal solution, which sometimes helps to improve interpretability of the results. Adding the additional $l_2$-regularizer helps to improve the performance when features are highly correlated~\cite{TWH15, ZH05}.}

\replaced{The convergence rates of iterative methods for solving~\eqref{eq:elastic} are typically governed by the condition number of \added[id=VM]{the Hessian matrix of the ridge loss, $\BB{C}+\reg_2\IM{}$, where} \replaced[id=VM]{$\BB{C}=\frac{1}{n}\Am^\top\Am$ is the  sample correlation matrix}{the  sample correlation matrix $\BB{C}=\frac{1}{n}\Am^\top\Am$}. Real-world data sets often have few  dominant features, while the other features are highly correlated with the stronger ones~\cite{TWH15,GOS16}.}{Consider the sample correlation matrix $\BB{C}=\frac{1}{n}\Am^\top\Am$, then the convergence rates of iterative methods for solving~\eqref{eq:elastic} are typically governed by the conditioning of $\BB{C}$.In practical data sets, we often have few strong dominant features, while other features are highly correlated with the stronger ones.} 
This translates to a \replaced{rapidly}{rapid} decaying spectrum of $\BB{C}$. In this paper, we  demonstrate how this property can be exploited to reduce the effect of ill-conditioning and to design faster algorithms for solving the elastic net regression problem~\eqref{eq:elastic}.

\subsection{Related work}
Over the past few years, there has been a great attention to \added{developing efficient optimization algorithms for minimizing composite objective functions} \begin{align}\label{eq:finite:sum}
	\min_{\xv \in \R^d} \FF{\xv} \triangleq \fF{}{\xv} +  \hF{\xv},
\end{align}
where 
$\fF{}{\xv}=\frac{1}{n}\sum_{i=1}^{n} \fF{i}{\xv} \label{eqn:finite-sum}$
is a finite sum of smooth \added{and convex} component functions $\fF{i}{\xv}$, and $\hF{\xv}$ is a possibly non-smooth convex regularizer. 
In machine learning applications, the function $f$ typically models the empirical data loss and the regularizer $h$ is used to promote desired properties of a solution. For example, the elastic net objective can fit to this form with 
$	f_i(\xv)
	=
		\frac{1}{2}\left(\av_i^\top\xv-\bv_i\right)^2
		+
		\reg_2 \Vnorm{\xv}_2^2,
$
and $h(\xv)=\reg_1\Vnorm{\xv}_1.$ 

\subsubsection{ First-Order methods}
Standard deterministic first-order methods for solving~\eqref{eq:finite:sum}, such as proximal gradient descent, enjoy linear convergence for strongly convex objective functions and are able to find an $\epsilon$-approximate solution in \added[id=VM]{time} $O\left(dn\kappa\log\frac{1}{\epsilon}\right)$\added{\deleted[id=VM]{iterations}, where $\kappa$ is the condition number {of $f$}. This\deleted[id=VM]{rate} \added[id=VM]{runtime} can}  be improved to $O\left(dn\sqrt{\kappa}\log\frac{1}{\epsilon}\right)$ if {it is combined with Nesterov acceleration}~\cite{BT09,Nes13}. 
\replaced{However, the}{The} main drawback of these methods is that they need to access the whole data set in every iteration, which {is too costly} in many machine learning tasks.

\added[id=VM]{For large-scale problems,  methods based on stochastic gradients have become the standard choice for solving~\eqref{eq:finite:sum}.} 
Many linearly convergent proximal methods\deleted{,} such as, \texttt{SAGA}~\cite{DBL14} \replaced{and}{or} \texttt{Prox-SVRG}~\cite{XZ14}, have been introduced and shown to outperform  standard first-order methods under certain regularity assumptions. These methods improve the time complexity to $O\left(d\left(n+\ktil\right)\log\frac{1}{\epsilon}\right)$, where $\ktil$ is a condition number satisfying $\ktil\geq\kappa$. 
When the component functions do not vary substantially in smoothness, $\ktil\approx\kappa$, \added{and} this complexity is far {better than} those of deterministic methods above. \replaced{By exploiting Nesterov}{Combining with Nesterov's} momentum in different ways (see, e.g.,~\cite{FGKS15b,LMH15,All16,Def16}), one can improve the complexity to $O\big(d\big(n+\sqrt{n\ktil}\big)\log\frac{1}{\epsilon}\big)$, which is also optimal for this class of problems~\cite{WS16}.

\subsubsection{Second-order methods}
Second-order methods are known to have superior performance compared to their first-order counterparts both in theory and practice, especially\deleted{,} when the problem at hand is high\added{ly} \replaced{nonlinear and/or ill-conditioned}{nonlinearity and/or ill-conditioning}. However, such methods often have very high computational cost per iteration\deleted[id=VM]{, due to the need of forming and solving a $d\times d$ linear system of equations.}.
Recently, there has been an intense effort to develop algorithms which use second-order information with a more reasonable computational burden~(see, e.g.,~\cite{ EM15,MNJ16,RM16,ABH17,XZ14,XYR16,BHNS16}~and references therein). Those methods use techniques such as random sketching, matrix sampling, and iterative estimation to construct an approximate Hessian matrix. Local and global \replaced{convergence guarantees}{convergences} have been derived  under various assumptions. Although many experimental results have shown excellent performance of those methods on {many} machine learning tasks, current second-order methods for finite-sum optimization tend to have much higher time-complexities than their first-order counterparts  (see~\cite{XYR16} for a detailed comparison).  

Apart from having high time complexities, none of the methods \replaced{cited}{mentioned} above \replaced{have any guarantees}{cannot work} in the composite setting since \replaced{their analyses {hinge} on differentiability of the objective function.}{ the objective functions are }\deleted{no longer}\deleted{differentiable.} \replaced{Instead, one has to rely on methods}{Few exceptions are methods} that build on \emph{proximal Newton} updates (see, e.g.,~\cite{LSS12, LHLS17, GS16, RK16}). \replaced[id=VM]{However, these methods still inherit the high update and storage costs of conventional second-order methods or require elaborate tuning of several parameters and stopping criteria depending on  a phase transition which occurs in the algorithm. }{The work in~\cite{LSS12} can be considered as an extension of the pure Newton method to the composite setting, hence having very costly updates. To reduce that cost, the work in~\cite{GS16} employs proximal quasi-Newton updates, while~\cite{RK16} proposes an incremental Newton method. However, the global convergence rates of these methods do not improve upon that of deterministic first-order methods, and they also need to store $d \times d$ matrices. In~\cite{LHLS17}, the update directions are approximated by an iterative routine. However, it requires to tune several parameters and stopping criteria depending on {\color{blue} a phase transition that occurs in the algorithm.} 
}

\subsubsection{Ridge regression}
For the smooth ridge regression problem, the authors in~\cite{GOS16} have developed a preconditioning method based on linear sketching which, when coupled with \texttt{SVRG}, yields a significant speed-up over stochastic first-order methods. This is a rare second-order method that has a comparable or even better time complexity than  stochastic first-order methods. More precisely, \replaced{it has a guaranteed}{they obtain the} running time \added{of}
$O(d(n + \kavg	) \log\frac{1}{\epsilon})$, where $\kavg$ is a new condition number that can be dramatically smaller than $\ktil$, especially when the spectrum of $\BB{C}$ decays rapidly. When $d \ll n$, the authors in \cite{WZ17} combine sub-sampled Newton methods with the mini-batch \texttt{SVRG} to obtain some further improvements.  
\subsection{Contributions}

{Recently, the work \cite{AS17} shows that under some mild
algorithmic assumptions, and if the dimension is sufficiently large, the iteration complexity of second-order methods for \emph{smooth} finite-sum problems composed of quadratics is no better than first-order methods. } 
Therefore, it is natural to ask whether one can develop a second-order method for solving the elastic net problem which \added{has improved practical performance but} still enjoys a strong \emph{worst-case} time complexity like the stochastic first-order methods do? 
It should be emphasized that due to the non-smooth objective, achieving this goal is much more challenging than for ridge regression. The preconditioning approach in~\cite{GOS16} is not applicable, and the current theoretical results for second-order methods are not likely to offer the desired  running time. 

In this paper, we provide a positive answer to this question. Our main contribution  is the design and analysis of 
a simple second-order method for the elastic net problem which has a strong theoretical time complexity \added{and superior practical performance}. 
{The convergence bound adapts to the problem structure and is governed by the spectrum and a \emph{statistical measure} of the data matrix. These quantities often yield significantly stronger time complexity guarantees {for practical datasets} than those of stochastic first-order methods (see Table~\ref{tab:comparison}).} 	
To achieve this, we first leverage recent advances in randomized low-rank approximation to generate a simple, one-shot approximation of the Hessian matrix. We then exploit the composite and finite-sum structure of the problem to develop a variant of the \texttt{ProxSVRG} method that builds upon\deleted{the} Nesterov's momentum acceleration and inexact computations of \emph{scaled proximal} operators, which may be of independent interest. We provide a simple convergence proof based on an explicit Lyapunov function, thus avoiding the use of sophisticated \emph{stochastic estimate sequences}. 


\begin{table}[t]
\caption{
		\footnotesize{Summary of different algorithms solving the elastic net problem. Here, $\kappa$ and $\ktil$ are conventional condition numbers satisfying $\kappa \leq \ktil$, while $\kavg$ is a new condition number defined w.r.t the $\He$-norm. When $\He$ is an approximate Hessian of the ridge loss, $\kavg$ is often much smaller than $\ktil$, especially on practical data sets.}
	}
	\label{tab:comparison}
	\begin{center}
		\begin{tabular}{c|c|c}
			\hline 
			Algorithm 	& Time complexity 		& 2nd-order	
			\bigstrut\\
			\hline\hline
			\texttt{PGD}
			& $O(dn\kappa\log\frac{1}{\epsilon})$
			& no
			\\[0.75ex]
			\hline			
			\texttt{FISTA}
			& $O\left(dn\sqrt{\kappa}\log\frac{1}{\epsilon}\right)$
			& no
			\\[0.75ex]
			\hline
			\texttt{ProxSVRG}
			& $O\left(d\left(n+\ktil\right)\log\frac{1}{\epsilon}\right)$
			& no
			\\[0.75ex]
			\hline
			\texttt{Katyusha}
			& $O\big(d(n+\sqrt{n\ktil})\log\frac{1}{\epsilon}\big)$
			& no
			\\[0.75ex]
			\hline
			\texttt{Ours}
			& $O\left(d\left(n+\kavg\right)\log\frac{1}{\epsilon}\right)$
			& yes
			\\[0.1ex]
			\hline
		\end{tabular}
	\end{center}
\end{table}

\section{Preliminaries and Notation}
Vectors are indicated by bold lower-case letters, and matrices are denoted by bold upper-case letters. We denote the dot product between $\xv$ and $\yv$ by  $\InP{\xv}{\yv}=\xv^\top\yv$, and the Euclidean norm of $\xv$ by $\Vnorm{\xv}_2=\sqrt{\InP{\xv}{\xv}}$. For a symmetric positive definite matrix $\He$, $\InP{\xv}{\yv}_{\He} =\xv^\top\He\yv$ is the $\He$-inner product of two vectors $\xv$ and $\yv$ and $\Mnorm{\xv}{\He}=\sqrt{\InP{\xv}{\xv}_{\He}}$ is the Mahalanobis norm of $\xv$. We denote by $\eig_i\left(\Am\right)$ the $i$th largest eigenvalue of $\Am$. Finally, $\eig_i$ denotes the $i$th largest eigenvalue of the correlation matrix $\BB{C}$.

In the paper, we shall frequently use the notions of strong convexity and smoothness in the $\He$-norm, introduced {in the next two assumptions}.
\begin{assumption}\label{asmp:1}
The function $\hF{\xv}$ is lower-semicontinous and convex and  $\mathrm{dom}\, h:= \{\xv \in \R^d \, | \, \hF{\xv} < \infty\}$, is closed. Each function $f_i$ is $L_i$-smooth w.r.t the $\He$-norm, i.e, there exists a positive constant $L_i$ such that
\begin{align*}
	\Mnorm{\grad{\fF{i}{\xv}}-\grad{\fF{i}{\yv}}}{\He^{-1}} \leq L_i \Mnorm{\xv-\yv}{\He}, \,\, \forall \xv, \yv \in \R^d.
\end{align*} 
\end{assumption}

Assumption~\ref{asmp:1} implies that $\nabla f$ is $L$-Lipschitz:
\begin{align*}
	\Mnorm{\grad{\fF{}{\xv}}-\grad{\fF{}{\yv}}}{\He^{-1}}
	\leq
	L
	\Mnorm{\xv-\yv}{\He}
\end{align*}
for some $L\leq\Lavg=\frac{1}{n}\sum_{i=1}^{n}L_i$. As a consequence, we have the following bound:
\begin{align*}
	\fF{}{\yv} \leq \fF{}{\xv} + \InP{\grad{\fF{}{\xv}}}{\yv-\xv} +  \frac{\Lavg}{2}\Mnorm{\yv-\xv}{\He}^2.
\end{align*}

\begin{assumption}\label{asmp:2}
The function $\fF{}{\xv}$ is $\mu$-strongly convex w.r.t the $\He$-norm, i.e, there exists a positive constant $\mu$ such that
\begin{align*}
	\fF{}{\lambda\xv+\left(1-\lambda\right)\yv}
	&\leq
		\lambda\fF{}{\xv}
		+
		\left(1+\lambda\right)\fF{}{\yv}
		-
		\frac{\mu\lambda\left(1-\lambda\right)}{2}
		\Mnorm{\xv-\yv}{\He}^2
\end{align*}
holds for all $x, y \in \R^n$ and $\lambda\in\left[0,1\right]$. 
\end{assumption}

Assumption~\ref{asmp:2} is equivalent to the requirement that
\begin{align*}
	\fF{}{\yv} \geq \fF{}{\xv} + \InP{\grad{\fF{}{\xv}}}{\yv-\xv} + \frac{\mu}{2}\Mnorm{\yv-\xv}{\He}^2, 
\end{align*} 
holds for all $\xv, \yv \in \R^n$. We will use both of these definitions of strong convexity in our proofs.

At the core of our method is the concept of \emph{scaled proximal mappings}, defined as follows:
\begin{definition}[Scaled Proximal Mapping]
For a  convex function $h$ and a symmetric positive definite matrix $\He$, the scaled proximal mapping of $h$ at $\xv$ is
\begin{align}\label{eq:def:prox}
	\Prox{h}{\BB{H}}{\yv}
	=
	\argmin_{\xv \in \R^d}\big\{\hF{\xv} + \frac{1}{2} \Mnorm{\xv-\yv}{\BB{H}}^2\big\}.
\end{align}	
\end{definition}
The scaled proximal mappings generalize the conventional ones:
\begin{align}\label{eq:def:l2_prox}
	\Prox{h}{}{\yv}
	=
	\argmin_{\xv \in \R^d}\big\{\hF{\xv} + \frac{1}{2} \Mnorm{\xv-\yv}{2}^2\big\}.
\end{align}	
 However, while many conventional prox-mappings admit analytical solutions, \replaced{this is almost never the case for scaled proximal mappings. This makes it hard to extend efficient first-order proximal methods to second-order ones.}{it almost never happens in the case of scaled proximal mappings making extending first-order proximal methods to second-order ones very challenging. }
{Fortunately}, scaled proximal mappings do share some key properties with the conventional ones. We collect a few of them in the following result:
\begin{property}[\cite{LSS12}]\label{lem:prox}
The following properties hold:
\begin{description}
\item  1.  $\Prox{h}{\He}{\xv}$ exists and is unique for $\xv \in \mathrm{dom}\, h$.
\item 2. Let $\partial h\left(\xv\right)$ be the subdifferential of $h$ at $\xv$, then 
\begin{align*}
	\He\left(\xv-\Prox{h}{\He}{\xv}\right)
	\in 
		\partial h\left(\Prox{h}{\He}{\xv}\right).
\end{align*}
\item  3. $\Prox{h}{\He}{\cdot}$ is non-expansive in the $\He$-norm:
\begin{align*}
	\Mnorm{\prox_{h}^{\He}(\xv)-\prox_{h}^{\He}(\yv)}{\He} 
	\leq
	\Mnorm{\xv-\yv}{\He}
	\,\, 
	\forall \xv, \yv \in \mathrm{dom}\, h.
\end{align*}
\end{description}
\end{property}

Finally, in our algorithm, it will be enough to solve~\eqref{eq:def:prox} approximately in the following sense:
\begin{definition}[Inexact subproblem solutions]
We say that $\xv^+ \in \R^d$ is an $\epsilon$-optimal solution to~\eqref{eq:def:prox} if
\begin{align}\label{def:eps:sol}
	&\hF{\xv^+} + \frac{1}{2\eta} \Mnorm{\xv^+-\yv}{\BB{H}}^2
	\leq
	\underset{\xv\in\R^d}{\min} 
		\big\{\hF{\xv} + \frac{1}{2\eta} \Mnorm{\xv-\yv}{\BB{H}}^2\big\}
	+
	\epsilon.
\end{align}

\end{definition}

\section{Building Block 1: Randomized Low-Rank Approximation}
The computational cost of many Newton\added{-}type methods is\deleted{often} dominated by the time required to compute the update direction\deleted{, namely,} $\dv=\BB{H}^{-1}\gv$ for some vector $\gv \in \R^d$ and approximate Hessian $\BB{H}$. A naive implementation using SVD would take $O\left(nd^2\right)$ flops,\deleted{and} \replaced{which is prohibitive for}{hence prohibited from} large-scale data sets. A natural way to reduce this cost is to use  truncated SVD.  However, standard deterministic methods  such as the power method and the Lanzcos method have run times that {scale inversely with} the gap between the eigenvalues of the input matrix. \replaced[id=VM]{This gap can be arbitrarily small for practical data sets, thereby preventing us from obtaining the desired time complexity.}{which \replaced{leads to poor performance on}{are often very small in} practical large data sets.} In contrast,\deleted{the} randomized sketching schemes usually admit gap-free run times~\cite{HMT11}. However, unlike other methods, the block Lanczos method, detailed in Algorithm~\ref{alg:BLZ}, admits both fast run times and  
strong guarantees on the errors between the true and the computed approximate singular vectors.  This property turns out to be critical for deriving bounds on the condition number\deleted{s} of the elastic net. 
\begin{proposition}[\cite{MM15}]
Assume that ${\BB{U}}_r$, ${\BB{\Sigma}}_r$, and ${\BB{V}}_r$ are matrices generated by Algorithm~\ref{alg:BLZ}. Let ${\Am}_r={\BB{U}}_r{\BB{\Sigma}}_r{\BB{V}}_r^\top=\sum_{i=1}^{r}{\sigma}_i{\BB{\uv}}_i{\BB{\vv}}_i^\top$ and let  $\Am=\sum_{i=1}^{d}\bar{\sigma}_i\bar{\uv}_i\bar{\vv}_i^\top$ be the SVD of $\Am$. Then, the following bounds hold with probability at least $9/10$:
\begin{align*}
	\Vert\Am-{\Am}_r\Vert_2
	&\leq
		\left(1+\epsilon^{\prime}\right)
		\bar{\sigma}_{k}	
	\\
	\big|{\BB{\uv}}_i^\top\Am\Am^\top{\BB{\uv}}_i-\bar{\uv}_i^\top\Am\Am^\top\bar{\uv}_i\big|
	&\leq		
		\epsilon^{\prime}\bar{\sigma}_{r+1}^2, \quad\forall i \in \{1,\ldots,r\}.
\end{align*} 
The total running time is \replaced{$O\left(ndr\log d(\epsilon^{\prime})^{-1/2}\right)$}{$O\left(\frac{ndr\log d}{\sqrt{\epsilon'}}\right)$}.
\end{proposition}

{Note that we only run Algorithm~1 once and $\epsilon'=1/2$ is sufficient in our work. Thus, the computational cost of this step is negligible, in theory and in practice.} 

\added[id=VM]{\subsection{Aproximating the Hessian}}
In this work, we consider {the following} approximate Hessian matrix {of the ridge loss}:
\begin{align}\label{eq:approx:hess}
	\He
	=
		{\BB{V}}_r
		\left({\BB{\Sigma}}_r^2+\reg_2\IM{}\right)
		{\BB{V}}_r^\top
		+
		\left(\sigma_r^2+\reg_2\right)
		\left(\IM{}-{\BB{V}}_r{\BB{V}}_r^\top\right).
\end{align}
\added[id=VM]{Here, the first term is a natural rank $r$ approximation of the true Hessian, while the second term is used to capture information in the subspace orthogonal to the column space of $\BB{V}_r$. The inverse of $\He$ in \eqref{eq:approx:hess} admits the explicit expression}
\begin{align*}
	\He^{-1}
	=
		{\BB{V}}_r
		\left({\BB{\Sigma}}_r^2+\reg_2\IM{}\right)^{-1}
		{\BB{V}}_r^\top
		+
		\frac{1}{\sigma_r^2+\reg_2}
		\left(\IM{}-{\BB{V}}_r{\BB{V}}_r^\top\right),
\end{align*}
so the evaluation of $\He^{-1}\xv$ has time complexity $O\left(rd\right)$.

\begin{algorithm}[!t]
\caption{Randomized Block Lanczos Method~\cite{MM15}}

\begin{algorithmic}[1]\label{alg:BLZ}
\REQUIRE Data matrix $\Am\in\R^{n\times d}$, target rank $r$, target precision $\epsilon'\in\left(0,1\right)$ \\[0.1cm]
\STATE Let $q=O(\log d/\sqrt{\epsilon'})$, and draw $\BB{\Pi} \sim \mathcal{N}_{d \times r}\left(\B{0},\IM{}\right)$   \\[0.1cm]
\STATE Compute $\BB{K}=\left[\Am \BB{\Pi} \,\, \left(\Am\Am^\top\right)\Am \BB{\Pi}\,\, \ldots\,\, \left(\Am\Am^\top\right)^q\Am \BB{\Pi}\right]$  \\[0.1cm]
\STATE Orthonormalize columns of $\BB{K}$ to obtain $\BB{Q}$\\[0.1cm]
\STATE Compute truncated $r$-SVD of $\BB{Q}^\top\Am$ as ${\BB{W}}_r{\BB{\Sigma}}_r{\BB{V}}_r^{\top}$
\STATE Compute ${\BB{U}}_r=\BB{Q}{\BB{W}}_r$
\ENSURE ${\BB{U}}_r$, ${\BB{\Sigma}}_r$, ${\BB{V}}_r$
\end{algorithmic}
\end{algorithm}

\subsection{Bounding the Condition Number\deleted[id=VM]{s}}
We now turn our attention to studying how the approximate Hessian affects the relevant condition number of the elastic net problem.  
We first introduce a condition number that usually determines the iteration complexity of stochastic first-order methods \added[id=VM]{under non-uniform sampling}.
\begin{definition}
The average condition number \added{of (\ref{eqn:finite-sum})} is
\begin{align*}
	\kavg=\frac{\Lavg}{\mu}=\frac{\frac{1}{n}\sum_{i=1}^{n}L_i}{\mu}.
\end{align*}
\end{definition}

For the elastic net problem~\eqref{eq:elastic}, the smooth part \replaced{of the objective is the ridge loss}{of its objective in can be written as}
\begin{align*}
		\frac{1}{n}
		\sum_{i=1}^{n}
		\underbrace{\frac{1}{2}\left(\av_i^\top\xv-\bv_i\right)^2
		+
		\reg_2 \Vnorm{\xv}_2^2}_{f_{i}(\xv)}.
\end{align*}
Since we define smoothness and strong convexity of $f_i(\xv)$ in the $\He$-norm, the relevant constants are 
\begin{align*}
	L_i&=\Vert\He^{-1}\big(\av_i\av_i^\top+\reg_2\IM{}\big)\Vert_2\\
	\mu&=\eig_d\big(\He^{-1/2}\left(\BB{C}+\reg_2\IM{}\right)\He^{-1/2}\big).
\end{align*}

\replaced{For comparison,}{For a comparison purpose} we  also \replaced{define}{use} the conventional condition number $\ktil$, which characterizes the smoothness and strong convexity of $f_i(\xv)$ in the Euclidean norm. In this case  $\tilde{\kappa}=\sum_i L_i /(n\mu)$, where
\begin{align*}
	L_i=\Vert\av_i\av_i^\top+\reg_2\IM{}\Vert^2_2 
	\quad\mbox{and} \quad
	\mu=\eig_d\left(\BB{C}+\reg_2\IM{}\right).
\end{align*}

\replaced[id=VM]{It will become apparent that $\kavg$ can be expressed in terms of a statistical measure of the ridge loss and that it may be significantly smaller than $\ktil$. 
We start by introducing a statistical measure that \replaced{has been}{were} widely used in the analysis of ridge regression (see, e.g., \cite{HKZ12} and the references therein).}{It will become apparent that both $\kmax$ and $\kavg$ can be expressed in terms of statistical measures of the ridge loss and \added{that they may be} significantly smaller than the conventional ones. 
We start by introducing \replaced{two}{the}\deleted{notions of} statistical measure\added{s} that were widely used in the analysis of ridge regression (see, e.g., \mbox{\cite{HKZ12,WZ17}} and \added{the} references therein).}
\begin{definition}[Effective Dimension]
For a positive constant $\lambda$, the effective dimension of $\BB{C}$ is defined as
\begin{align*}
	\sdim_{\lambda}=\sum_{i=1}^{d}\frac{\eig_i}{\eig_i+\lambda}.
\end{align*}
\deleted{
The statistical leverage at $\lambda$ of the data matrix $\Am$ is said to be bounded by $\lvr_{\lambda}$  if }
\end{definition}
\replaced[id=VM]{The effective dimension generalizes the ordinary dimension and satisfies  $\sdim_{\lambda}\leq d$ with equality if and only if $\lambda=0$.} {Clearly, $\sdim_{\lambda}\leq d$ with equality if and only if $\lambda=0$.} It is typical that when $\BB{C}$ has a \replaced{rapidly}{fast} decaying spectrum, most of the $\lambda_i$'s are dominated by $\lambda$, and hence $\sdim_\lambda$ can be much smaller than $d$. \deleted[id=VM]{For the statistical leverage, under some standard assumptions, it was shown in \mbox{\cite{HKZ12}} that $\lvr_{\lambda}$ only grows logarithmically with the dimensions of $\Am$. }

The following lemma\deleted{provides} bounds\deleted{for} the eigenvalues of \replaced[id=VM]{the matrix $\He^{-1/2}\left(\BB{C}+\reg_2\IM{}\right)\He^{-1/2}$, which can be seen as the effective Hessian matrix}{the effective Hessian matrix}.
\begin{lemma}[\!\!{\cite{GOS16}}]\label{lem:eig:bound}
Invoking Algorithm~\ref{alg:BLZ} with \replaced{data matrix}{parameters} $\frac{1}{\sqrt{n}}\Am$, \added{target rank} $r$, and \added{target precision} $\epsilon'=1/2$,\deleted{then} it holds with probability at least $9/10$ that
\begin{align*}
	\eig_1\left(\He^{-1/2}\left(\BB{C}+\reg_2\IM{}\right)\He^{-1/2}\right)
	\leq 17
	\\
	\frac{\reg_2}{19\left(\eig_r +\reg_2\right)}
	 \leq
	\eig_d\left(\He^{-1/2}\left(\BB{C}+\reg_2\IM{}\right)\He^{-1/2}\right)
	\leq 2.
\end{align*}
\end{lemma}

\replaced[id=VM]{Equipped with Lemma~\ref{lem:eig:bound}, we can now connect $\kavg$ with  $\sdim_{\lambda}$ using the following result.}{Equipped with Lemma~\ref{lem:eig:bound}, we can now connect $\kmax$ and $\kavg$ with the statistical measures $\sdim_{\lambda}$ and $\lvr_{\lambda}$ \replaced{using}{via} the following result.}
\begin{theorem}\label{thrm:bounding:cond}
\replaced{With probability at least $9/10$, the following bound holds up to a multiplicative constant:}{The following bounds hold up to multiplicative constants with probability at least $9/10$:}
\begin{align*}
	\kavg \leq \min\left(\frac{\sdim_{\reg_2}}{\reg_2},\frac{r\eig_r+\sum_{i>r}\eig_i}{\reg_2}+d\right).
\end{align*}	
\begin{proof}
See Appendix~\ref{apdx:thrm:bounding:cond}.
\end{proof}
\end{theorem}
 Since $\tilde{\kappa}=(\sum_i \lambda_i + d\gamma_2)/\gamma_2$,\deleted[id=VM]{As a consequence, \replaced{the maximum condition number is reduced by a factor of order}{when compared to the case $\He=\IM{}$, the gain in the maximum condition number is of order} $\frac{\max_i\Vnorm{\av_i}_2^2}{\lvr_{\reg_2}^2\sdim_{\reg_2}},$ \replaced{when operating in the $\BB{H}$-norm. This reduction}{which} can be significant since the value of $\sdim_{\lambda}\lvr_{\lambda}^2$ is small for \replaced{when the spectrum of $\BB{C}$ decays rapidly}{for fast decaying spectrum}, while $\max_i\Vnorm{\av_i}^2$ can grow unbounded with the dimension\replaced{. For}{, for} instance, when the samples $\av_i$ were draw from a sub-Gaussian distribution, $\max_i\Vnorm{\av_i}^2=O\left(d\right)$ \mbox{\cite{Ver18}}.}
%
$\kavg$ is reduced by a factor 
\begin{align*}
		\frac{
			\sum_{ i\leq r}\eig_i+\sum_{i>r}\eig_i
		}{
			r\eig_r+\sum_{i>r}\eig_i
		},
\end{align*}
 compared to $\tilde{\kappa}$.  
\replaced{If}{Therefore, if} the spectrum of $\BB{C}$ decays rapidly, then  the terms $\sum_{i>r} \lambda_i$ are negligible and the ratio is approximately $\sum_{i\leq r}\lambda_i / (r\lambda_r)$. If the first eigenvalues are much larger than $\lambda_r$, this ratio will be large. For example, for the australian data set~\cite{CL11}, this ratio can be as large as $1.34\times 10^4$ and $1.6\times 10^5$ for $r=3$ and $r=4$, respectively.
This indicates that it is possible to improve the iteration complexity of stochastic first-order methods if one can \replaced{capitalize on}{realize} the notions of strong convexity and smoothness w.r.t the $\He$-norm \replaced{in the}{by an} optimization algorithm. Of course, this is only meaningful if there is an efficient way to inject curvature information into the optimization process without \added{significantly increasing} the computational cost\deleted{ of stochastic first-order methods}. In the smooth case, i.e., $\reg_1=0$, this task can be done by a preconditioning step~\cite{GOS16}. However, this approach is not applicable for \added{the} elastic net, and we need to make use of another building block.

\section{Building Block 2: Inexact Accelerated Scaled Proximal SVRG}
In this section, we introduce an inexact scaled accelerated \texttt{ProxSVRG} \added{algorithm} for solving the generic finite-sum minimization problem in~\eqref{eq:finite:sum}. We then characterize the convergence rate of the proposed algorithm.
\subsection{Description of the Algorithm}
To motivate our algorithm, we first recall the \texttt{ProxSVRG} method \replaced{from}{in} \cite{XZ14}:
For the $s$th outer iteration with the \replaced{corresponding outer}{current} iterate $\tilde{\xv}_s$, \replaced{let}{define} $\xv_0=\tilde{\xv}_s$\replaced{ and}{, then} for $k=0,2,\ldots, T-1$\deleted{,} do
\begin{align}
	&\vv_k =  \left(\grad\fF{i_k}{\xv_k}-\grad\fF{i_k}{\tilde{\xv}_s}\right)/\left(np_{i_k}\right)+\grad\fF{}{\tilde{\xv}_s}
	\\\label{eq:prox:step:svrg}
	&\xv_{k+1} = 	\Prox{\eta h}{}{\xv_k - \eta\vv_k},
\end{align}
where $i_k$ is \replaced{drawn randomly}{randomly drawn} from $\{1,\ldots,n\}$ with\deleted{the} probability \replaced{$p_{i_k}=L_{i_k}/(n\Lavg)$}{$p_{i_k}=\frac{L_{i_k}}{n\Lavg}$}. 
Since we are provided with an approximate Hessian matrix $\He$, it is natural to use the following update:
\begin{align}\label{eq:prox:step}
	\xv_{k+1} = 	\Prox{\eta h}{\He}{\xv_k - \eta\He^{-1}\vv_k},
\end{align}
which can be seen as a proximal Newton step with the full gradient vector replaced by $\vv_k$. Note that when $h(\cdot)$ is the $\ell_1$-penalty, \texttt{ProxSVRG} can evaluate~\eqref{eq:prox:step:svrg} 
in time $O(d)$, while evaluating~\eqref{eq:prox:step} amounts to solving an optimization problem. It is thus is critical to keep the number of such evaluations small, which then translates into making a sufficient progress at each iteration. A natural way to achieve this goal is to reduce the variance of the noisy gradient $\vv_k$.
This suggests to use large mini-batches, i.e., instead of using a single component function $f_{i_k}$, we use multiple ones to form:
\begin{align*}
	\vv_k 
	=
		\frac{1}{b}
		\sum_{i_k\in\mathcal{B}_k}
		\big(\grad\fF{i_k}{\xv_k}-\grad\fF{i_k}{\tilde{\xv}_s}\big)/\left(np_{i_k}\right)+\grad\fF{}{\tilde{\xv}_s},
\end{align*}
where $\mathcal{B}_k \subset \{1,\ldots,n\}$ is a set of indices with cardinality $\norm{\mathcal{B}_k}=b$. 
It is easy to verify that $\vv_k$ is an unbiased estimate of $\grad\fF{}{\xv_k}$. \added[id=VM]{Notice that naively increasing the batch size makes the algorithm increasingly similar to its deterministic counterpart, hence inheriting a high-time complexity.} This makes it hard to retain the runtime of \texttt{ProxSVRG} in the presence of 2nd-order information.
 
\begin{algorithm}[!t]
\caption{Inexact Accelerated Scaled Proximal SVRG}
 \begin{algorithmic}[1]\label{alg:1}
\REQUIRE $\xTd_0$, $\{\Bs_k\}_{k=0}^T$, $\eta$, $\tau$\\[0.1cm]
\FOR{$s=0,1,\ldots,S$}
\STATE  $\grad{\fF{}{\xTd_{s}}} \gets \frac{1}{n}\sum_{1}^{n}\grad{\fF{i}{\xTd_{s}}}$ 
\STATE $\xv_0 \gets \zv_0 \gets \xTd_s$ 
\FOR{$k=0,1,\ldots, T-1$}
\STATE $\yv_k\gets \frac{1}{1+\tau} \xv_k + \frac{\tau}{1+\tau} \zv_k\ $
\STATE $\vv_k \gets \grad{\fF{\Bs_k}{\yv_k}} - \grad{\fF{\Bs_k}{\xTd_{s}}} + \grad{\fF{}{\xTd_{s}}} $
\STATE  $\xv_{k+1} \approx \Prox{\eta h}{\He}{\yv_k- \eta \He^{-1}\vv_k}$
\STATE  $\gv_{k+1} \gets \frac{1}{\eta}\left(\yv_k - \xv_{k+1}\right) $
\STATE $\zv_{k+1}\gets \zv_k + \tau\left(\yv_k-\zv_k\right) - \frac{\tau}{\mu}\gv_{k+1}$
\ENDFOR
\STATE $\xTd_{s+1}\gets\xv_{T}$
\ENDFOR
\ENSURE $\xTd_{S}$
\end{algorithmic}
\end{algorithm} 
 
In the absence of second-order information \replaced{and under the assumption that}{as well as when} the proximal step is computed exactly, the work~\cite{Nit14} introduced a method called \texttt{AccProxSVRG} that enjoys the same time complexity \replaced{as}{of} \texttt{ProxSVRG} but  allows for much larger mini-batch sizes.
In fact, it can tolerate a mini-batch of size $O\big(\sqrt{\ktil}\big)$ \replaced{thanks to the use of Nesterov}{, and the key ingredient is the use of Nesterov's} momentum.  This indicates that an appropriate use of Nesterov momentum in our algorithm could allow for the larger mini-batches required to balance the computational cost of using scaled proximal mappings. The improved iteration complexity of the scaled proximal mappings will then give an overall acceleration in terms of wall-clock time. 
As discussed in~\cite{All16}, the momentum mechanism in \texttt{AccProxSVRG} fails to accelerate  \texttt{ProxSVRG} unless $\ktil \geq n^2$.  In contrast, as we will see, our algorithm will be able to accelerate the convergence also in these scenarios.
\added[id=VM]{In summary, our algorithm } is designed to run in an inner-outer fashion as \texttt{ProxSVRG} with large mini-batch sizes and Nesterov momentum to compensate for the increased computational cost of subproblems.
The overall procedure is summarized in Algorithm~\ref{alg:1}.

\subsection{Convergence Argument}

In this subsection, we will show that as long as the errors in evaluating the scaled proximal mappings are controlled in an appropriate way, the iterates generated by the outer loop\deleted{s} of Algorithm~\ref{alg:1} converge linearly \added[id=VM]{in expectation} to the optimal solution.
Recall that \replaced{in}{at} Step~7 of Algorithm~\ref{alg:1}, we want to find an $\epsilon_k$-optimal solution in the sense of~\eqref{def:eps:sol} to the following problem:
\begin{align}\label{eq:subproblem}
	\underset{\xv \in \R^d}{\mbox{minimize}}\,\frac{1}{2\eta} \Mnorm{\xv-\yv_k + \eta\He^{-1} \vv_k}{\He}^2 + \hF{\xv}.
\end{align}
\replaced{The}{Based on that, the} next lemma quantifies the progress made by one inner iteration of the algorithm.\deleted[id=VM]{\mbox{\footnote{Note that for the rest of the paper, we only consider the non\added{-}uniform sampling scheme, and all the subsequent results follow trivially for uniform sampling by setting $p_{i}=1/n$ and replacing $\Lavg$ by $\Lmax$ everywhere.}}}
Our proof builds on \replaced{a}{the} Lyapunov argument using \replaced{a Lyapunov function on}{the Lyapunov functions of} the form:
\begin{align}\label{eq:Lyp:def}
	\Lyp_k= \FF{\xv_k}-\FF{\xv^\star} + \frac{\mu}{2}\Mnorm{\zv_k-\xv^\star}{\He}^2.
\end{align}

\begin{lemma}\label{keylem}
Let Assumptions~\ref{asmp:1}--\ref{asmp:2} hold and let $\xv^\star = \argmin_{\xv} \FF{\xv}$\replaced{,}{. 
Let} $\eta = 1/\Lavg$ and\deleted{let} $\tau=\sqrt{\mu/2\Lavg}$. If the mini-batch size is chosen such that $b\geq 60\sqrt{\Lavg/\mu}$, then for any $k \in \{0,\ldots,T-1\}$, there exists a vector $\BB{\xi}_k \in \R^d$ such that $\Mnorm{\BB{\xi}_k}{\He^{-1}} \leq \sqrt{2\eta\epsilon_k}$ and
\begin{align}\label{eq:keylem}
	\E\Lyp_{k+1}
	&\leq	
		\left(1-\tau\right)\E \Lyp_k		
		+
		\tau \Lavg \E\InP{\BB{\xi}_k}{\xv^\star-\zv_k}
		+
		5\epsilon_k
	\nonumber\\
	&\hspace{0.5cm}
		+
		\frac{\tau}{5}
		\EX{
			\FF{\xv_k} - \FF{\xv^\star}
			+
			\FF{\xTd_s} - \FF{\xv^\star}
		}.
\end{align}
\begin{proof}
See Appendix~\ref{apdx:keylem}.
\end{proof}
\end{lemma}

\begin{remark}
Our proof \replaced{is direct and based}{bases directly} on\deleted{the rather} natural Lyapunov functions\replaced{, thereby avoiding}{ and avoids} the use of stochastic estimate subsequences as in~\cite{Nit14} which is already very complicated even when the subproblems are solved exactly and $\He=\IM{}$. 
We stress that the result in Lemma~\ref{keylem} also holds for smaller mini-batch sizes, namely $b\in\big\{1,\ldots,O(\sqrt{\Lavg/\mu})\big\}$, provided that the step size $\eta$ is reduced accordingly. In favor of a simple proof, we \replaced{only report}{opt to report only} the large mini-batch result here.
\end{remark}

\replaced{Equipped}{Provided} with Lemma~\ref{keylem}, we can now characterize the progress made by one outer  iteration of Algorithm~\ref{alg:1}.
\begin{theorem}\label{thrm:main}
Let Assumptions~\ref{asmp:1}--\ref{asmp:2} hold. Suppose that the parameters $\eta$, $b$, and $\tau$ are chosen according to Lemma~\ref{keylem}\replaced{ and define}{. Let} $\rho=9\tau/10$\replaced{. Then,}{ then,} if the errors in solving the subproblems satisfy  
\begin{align*}
	\epsilon_k 
	\leq 
		\left(1-\rho\right)^k\Lyp_0
\end{align*}
for all $k \in \{0,\ldots,T-1\}$\deleted{,} and 
$	T\geq (4\log c)/3\rho,$
where $c$ is a universal constant\replaced{, then for every $s\in \N_+$,}{.
Then, for any $s \in \N_+$, it holds that}
\begin{align*}
	\EX{
		\FF{\xTd_s}-\FF{\xv^\star}
	}
	\leq
	\frac{2}{3}
	\EX{\FF{\xTd_{s-1}}-\FF{\xv^\star}}.
\end{align*}
\begin{proof}
See Appendix~\ref{apdx:thrm:main}.
\end{proof}
\end{theorem}

\begin{remark}
The theorem indicates that if the errors in solving the subproblems are controlled appropriately, the outer iterates generated by Algorithm~\ref{alg:1} converge linearly \added{in expectation} to the optimal solution\deleted{ in expectation}. Since $\Lyp_0$ depends on $\xv^\star$, \replaced{it is difficult to provide a general closed-form expression for the}{in general, we cannot evaluate the} target precisions $\epsilon_k$. \deleted{However, for the elastic net problem, $\Lyp_0$ is  bounded by $2\FF{\xv_0}$, and hence, one can use this bound in place of $\Lyp_0$.}However, we will show below that with a certain \added{policy for selecting the} initial point, it is sufficient to run the solver a constant number of iterations independently of $\xv^\star$. {We stress that the results in this section are valid for minimizing general convex composite functions~\eqref{eq:finite:sum} and not limited to the elastic net problem.}
\end{remark}

\section{Warm-Start}
The overall time complexity of Algorithm~\ref{alg:1}\deleted{strongly} depends \added{strongly} on \replaced{our}{the} ability to  solve problem~\eqref{eq:subproblem} in a reasonable computational time. If one naively starts the solver at a random point, it may take many iterations to meet the target precision. Thus, it is necessary to have a well-designed warm-start procedure for initializing the solver. Intuitively, the current iterate $\xv_k$ can be a reasonable starting point since the next iterate $\xv_{k+1}$ should not be too far away from  $\xv_k$. 
However, in order to achieve a strong theoretical running time, we use a \replaced{rather}{relatively} different scheme inspired by~\cite{LMH15}. Let us first define the vector $\uv_k=\yv_k - \eta\He^{-1} \vv_k$ for $k\in\{0,1,\ldots,T-1\}$ and the function 
\begin{align*}
	p\left(\zv,\uv \right)=\hF{\zv}+\frac{1}{2\eta} \Mnorm{\zv-\uv}{\He}^2.
\end{align*}
Then, the $k$th subproblem seeks for $\xv_{k+1}$ such that
\begin{align}\label{eq:eps:req}
	p\left(\xv_{k+1},\uv_k \right)-p\left(\xv_{k+1}^\star,\uv_k \right)\leq \epsilon_k,
\end{align}
where $\xv_{k+1}^\star$ is the exact solution.
We consider \replaced{the}{an} initialization policy 
\begin{align}\label{eq:def:z0:main}
	\zv_0=\prox_{\gamma h}\big(\xv_k-\frac{\gamma}{\eta}\He\left(\xv_k-\uv_{k-1}\right)\big),
\end{align}
which can be seen as one step of the proximal gradient method applied to $p\left(\zv,\uv_{k-1} \right)$ starting at the current $\xv_k$ with\deleted{the} step size $\gamma$. 

The following proposition characterizes the \replaced{difference in objective realized by}{distance from} $\zv_0$ \replaced{and}{to} $\xv_{k+1}^\star$\replaced{.}{ in objective.}
\begin{proposition}\label{prop:ws}
 Let $\zv_0$ be defined by (\ref{eq:def:z0:main}) with $\gamma=\eta/\eig_1\left(\He\right)$. Let $\ksub=\eig_1\left(\He\right)/\eig_r\left(\He\right)$ be the condition number of the subproblems. Assume that the errors in solving the subproblems satisfy $\epsilon_k\leq\left(1-\rho \right)^k\Lyp_0$ for all $k\in\{0,1,\ldots,T-1\}$. Then,  
\begin{align*}
	p\left(\zv_0,\uv_k \right)-p\left(\xv_{k+1}^\star,\uv_k \right)
	\leq
		\frac{\ksub}{1-\rho}\epsilon_k.
\end{align*}	
\begin{proof}
See Appendix~\ref{apdx:prop:ws}.
\end{proof}
\end{proposition}

The proposition, together with~\eqref{eq:eps:req}, implies that it suffices to find $\xv_{k+1}$ such that 
\begin{align}\label{eq:restart:suff:cond}	
		&p\left(\xv_{k+1},\uv_k \right)-	p\left(\xv_{k+1}^\star,\uv_k \right)
	\leq 
		\frac{1-\rho}{\ksub}
	\left(	p\left(\zv_0,\uv_k \right)-	p\left(\xv_{k+1}^\star,\uv_k \right)\right).
\end{align}
This is significant since one only needs to reduce the residual error by a constant factor independent of the target precision. \added[id=VM]{Note also that the condition number  $\ksub$ is much smaller than $\kappa \approx \eig_1(\He)/(\eig_d+\gamma_2)$, and computing the gradient of the smooth part of $p\left(\zv,\uv \right)$ only takes time $O\left(rd\right)$ instead of $O\left(nd\right)$ as in the original problem. Those properties imply that the subproblems can be solved efficiently by iterative methods, where only a  small (and known) constant number of iterations is needed. The next section develops the final details of convergence proof.}

\section{Global Time Complexity}


We start with  the time complexity of Algorithm~\ref{alg:1}.
Let $\mathcal{T}\left(\alpha\right)$ be the number of gradient evaluations that a subproblem solver takes to reduce the residual error by a factor $\alpha$. Then, by Proposition~\ref{prop:ws}, one can find an $\epsilon_k$-optimal solution to the $k$th subproblem by at most  $\mathcal{T}\left(\ksub/\left(1-\rho\right)\right)$ gradient evaluations, where one gradient evaluation is equivalent to $O\left(d\right)$ flops.
Consider the same setting of Theorem~\ref{thrm:main} and suppose that the subproblems are initialized by~\eqref{eq:def:z0:main}. Then, the time complexity of Algorithm~\ref{alg:1} is given by: 
\begin{align}\label{eq:saproxsvrg:complx}
	O\big(d\big(n + \kavg + \sqrt{\kavg}\, \mathcal{T}\big(\frac{\ksub}{1-\rho}\big)\big)\log\frac{1}{\epsilon}\big),
\end{align}  
where the first summand is due to the full gradient evaluation at each outer loop; the second one comes from the fact that one needs  $O(\sqrt{\kavg})$ inner iterations, each of which uses a mini-batch of size $O(\sqrt{\kavg})$; and the third one is the result of $O(\sqrt{\kavg})$ inner iterations, each of which solves a subproblem that needs 
 $\mathcal{T}\left(\ksub/\left(1-\rho\right)\right)$ gradient evaluations.
We can now put things together and state our main result.
\begin{proposition}\label{prop:glob:runtime}
Suppose that the approximate Hessian matrix $\He$ is given by~\eqref{eq:approx:hess}\replaced{ and that Algorithm~\ref{alg:1} is invoked with}{, then by invoking Algorithm~\ref{alg:1} with} $\fF{}{\xv}=\frac{1}{2n}\Vnorm{\Am\xv-\bv}_2^2+\frac{\reg_2}{2}\Vnorm{\xv}_2^2$ and $\hF{\xv}=\Vnorm{\xv}_1$. Assume further that the subproblems are solved by the accelerated proximal gradient descent method~\cite{Nes13,BT09}. {Our method can find an $\epsilon$-optimal solution in time}
\begin{align*}
		O\big(d\left(n + \kavg \right)\log\frac{1}{\epsilon}\big).
\end{align*}  
\begin{proof}
\replaced[id=VM]{The task reduces to evaluating the term $\mathcal{T}\left(\ksub/\left(1-\rho\right)\right)$ in~\eqref{eq:saproxsvrg:complx}. Recall that the iteration complexity of the accelerated proximal gradient descent method for minimizing the function $F(x)=f(x)+h(x)$, where $f$ is a smooth and strongly convex function and $h$  is a possibly non-smooth convex regularizer, initiliazed at $\xv_0$ is given by $\sqrt{\kappa}\log\frac{F(\xv_0)-F(\xv^\star)}{\epsilon}$, where $\kappa$ is the condition number. By invoking the about result with $F(\xv)=p\left(\xv,\uv_k \right)$, $\xv^\star=\xv_{k+1}^\star$, $\xv_0=\zv_0$, $\kappa=\ksub$, and $\epsilon$ is the right-hand side of~\eqref{eq:restart:suff:cond}, it follows that the number of iterations for each subproblem can be bounded by}
{The result follows directly from~\eqref{eq:saproxsvrg:complx} and the fact that the number of iterations for each subproblem can be bounded by}
\begin{align}\label{eq:sub:complexity}
	O\big(\sqrt{\ksub}\log\frac{\ksub}{1-\rho}\big).
\end{align}
In addition, each iteration takes time $O\left(rd\right)$ to compute the gradient implying the  time complexity
\begin{align*}
	O\big(
		d
		\big(
			n 
			+ 
			\kavg 
			+ r\sqrt{\ksub}\sqrt{\kavg}\log\frac{\ksub}{1-\rho}
		\big)
		\log\frac{1}{\epsilon}
	\big)
\end{align*}
Finally, since $\kavg \gg  \ksub$, the proof is complete.
\end{proof}
\end{proposition}
We can easily recognize that this time complexity has the same form \replaced{as the}{of} stochastic first-order methods \added{discussed} in Section~1.2.1 with the condition number $\ktil$ replaced by $\kavg$. It has been shown in Theorem~\ref{thrm:bounding:cond} that $\kavg$ can be much smaller than $\ktil$, especially, when $\BB{C}$ has a \replaced{rapidly}{fast} decaying spectrum. {We emphasize that the expression in~\eqref{eq:sub:complexity} is available for free to us after having approximated the Hessian matrix. Hence \emph{no tuning} is required to set $\epsilon_k$ for solving the subproblems.} 

\begin{figure}[t!]
	\centering 
	{\includegraphics[width=.8\textwidth]{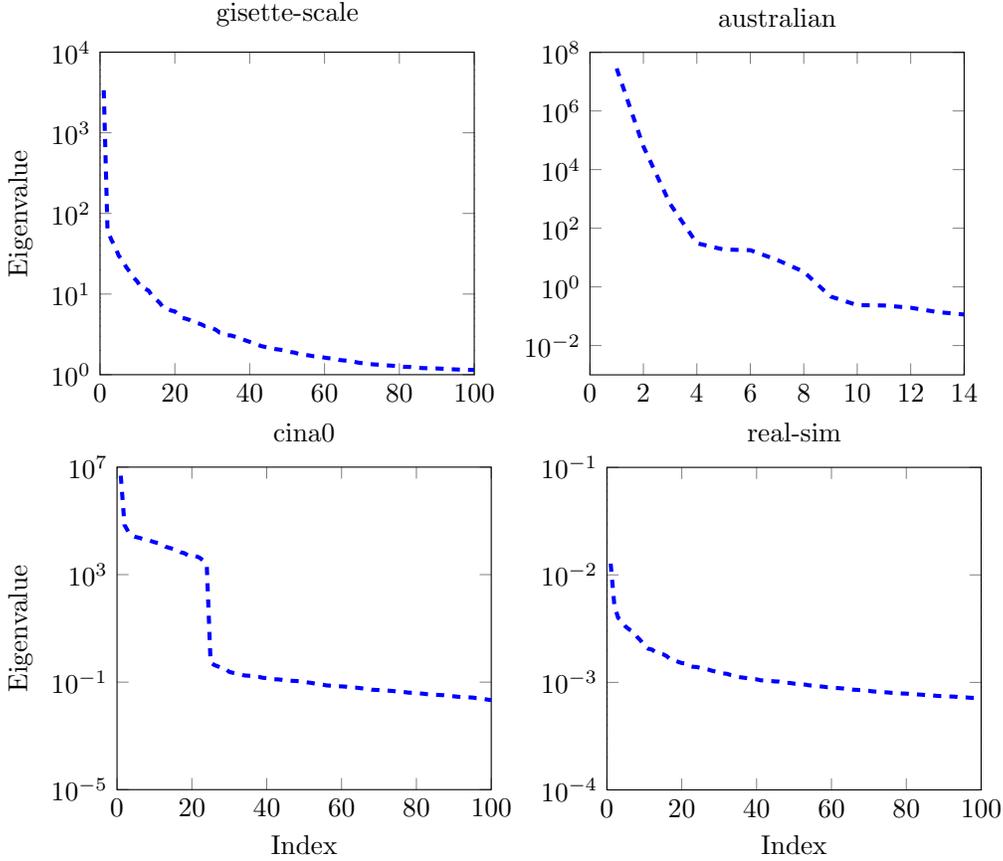}}
	\label{fig:spectra} 
	\caption{Spectrum of the empirical correlation matrix $\B{C}$ for different data sets.}
\end{figure}

\section{Experimental Results}\label{sec:expr}

In this section, we perform numerical experiments to verify the efficacy of the proposed method on real world data sets~\cite{CL11,GAC08}\deleted[id=VM]{ with various sizes and \replaced{condition numbers}{conditioning}. 
We fix $\reg_1$ as in Table~\ref{table:data} and vary $\reg_2$ from $\{\reg_1, 0.1\reg_1\}$ to influence the conditioning of the \replaced{ridge loss}{ problem}}.
We compare our method with several well-known first-order\deleted[id=VM]{ \added{mini-batch}} methods\deleted{ in mini-batch settings}:  \texttt{FISTA}~\cite{BT09} with optimal step-size;
 \texttt{Prox-SVRG}~\cite{XZ14} with epoch length $2n/b$ as suggested by the authors;
\texttt{Katyusha1}~\cite{All16} with epoch length $2n/b$, Katyusha momentum $\tau_2=0.5/b$ as suggested by the author;
 and our method with epoch length $2n/b$.
Since \texttt{Katyusha1} can use a mini-batch of size $\sqrt{n}$ without slowing down the convergence, we set $b=\sqrt{n}$ for all methods except \texttt{FISTA}. Finally, to make a fair comparison, for each algorithm above, we tune only the step size, from the set $\eta\times\{10^k, 2\times10^k, 5\times10^k| k\in\{0, \pm 1,\pm 2\}\}$, where $\eta$ is the theoretical step size, and report the one having smallest objective value. Other hyper-parameters are set to their theory-predicted values. For \texttt{Katyusha1}, we also compute the largest eigenvalue of the Hessian matrix in order to set its parameter $\tau_1$. All methods are initialized at $\BB{0}$ \replaced{and run for up to 100 epochs (passes through the full data set).}{. Number of epochs or data passes is up to  100.} For the subproblems in Algorithm~\ref{alg:1}, we just simply run \texttt{FISTA} with $\sqrt{\ksub}\log\ksub$ iterations as discussed in the previous section, without any further tunning steps. 
The value of $r$ is chosen as a small fraction of $d$ so that the preprocessing time of Algorithm~\ref{alg:BLZ} is negligible. Note that the available spectrum of $\BB{C}$ after running Algorithm~\ref{alg:BLZ} also provides an insightful way to choose $r$.
\begin{table}[h!]
    \caption{
                \footnotesize{Brief summary of data sets and parameters}
    }
    \label{table:data}
\begin{center}
    \begin{tabular}{c|c|c|c|c}
        \hline \hline
        Data set 	&  $d$ & $n$ & $\reg_1$ & $r$ 
        \bigstrut\\
        \hline
        gisette-scale	&	5,000 		&   6,000	&$10^{-3}$	&40			\\[0.3ex]
        australian  &	14	 		& 	690		&$10^{-3}$	&5		\\[0.3ex]
        cina0     	&	132 		&  16,033	&$10^{-4}$	&20			\\[0.3ex]
        realsim		&	20,958 		&  72,309	&$10^{-4}$	&50		\\[0.3ex]
        \hline			
    \end{tabular}
\end{center}
\end{table}

\begin{figure}[t!]
	\centering 
	{\includegraphics[width=.8\textwidth]{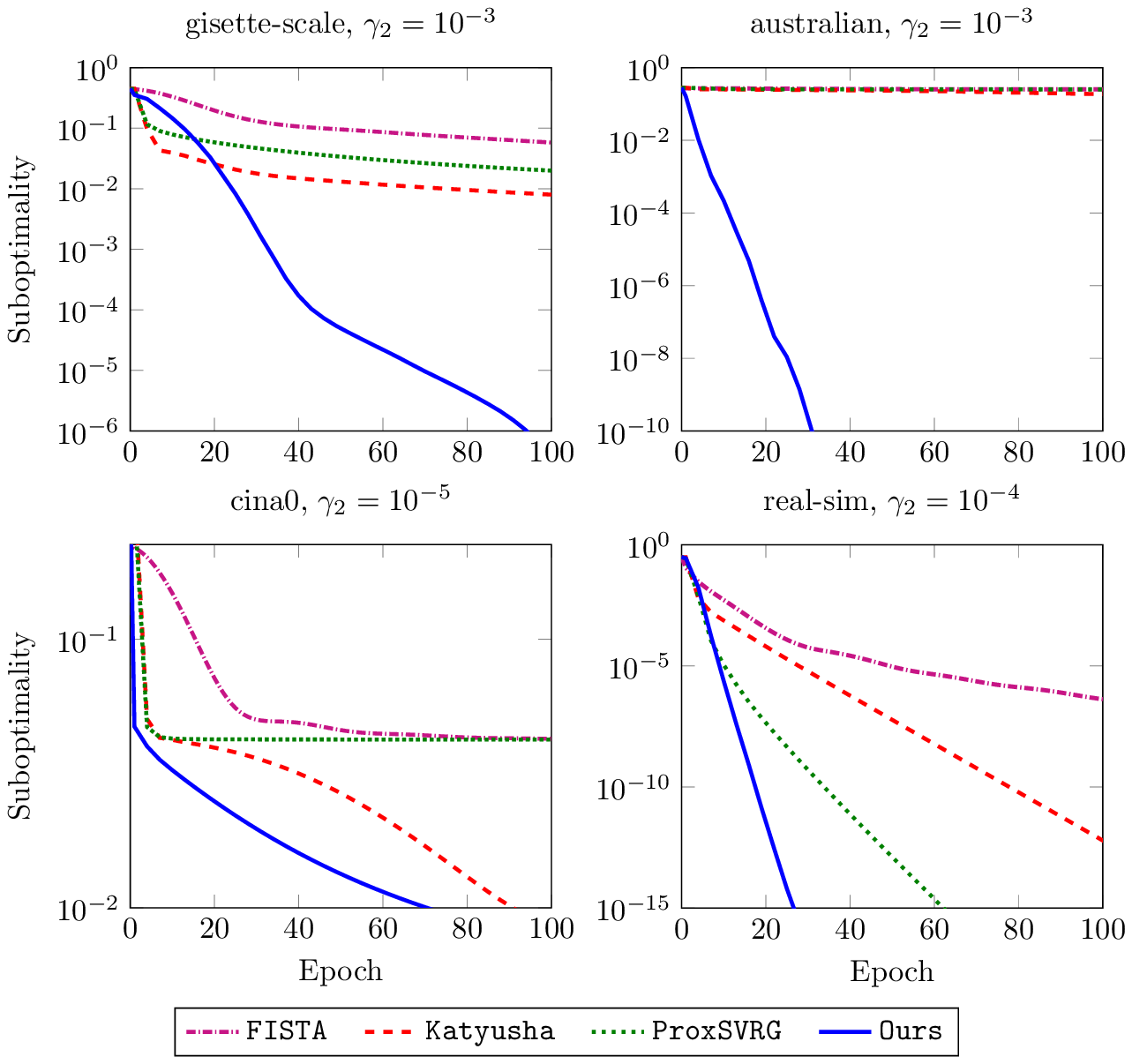}}
	\caption{Suboptimality versus the number of epochs for different algorithms solving the elastic net problem.}
	\label{fig:subopt} 
\end{figure}

\replaced[id=VM]{Figure~\ref{fig:subopt} shows the suboptimality in objective versus the number of epochs\deleted{(passes over the full data set)} for different algorithms solving the elastic net problem. We can see that our method systematically outperforms the others in all setting\added{s}, and that there is a clear correspondence between the spectrum of $\BB{C}$ in Fig.~\ref{fig:spectra} and the potential speed-up.   Notably, for the australian dataset, all the first-order methods make almost no progress in the first 100 epochs, while our method can find a high\added{-}accuracy solution within tens of epochs, demonstrating a great benefit of second-order information. On the other hand, for a well-conditioned data set that does not exhibit high curvature such as real-sim, \texttt{ProxSVRG} \replaced{is}{performs quite}  comparable to our method and even outperforms \texttt{Katyusha1}. This\deleted{highly} agrees with the theoretical time complexities summarized in Table~\ref{tab:comparison}. Finally, we demonstrate a hard instance for low-rank approximation methods via the cina0 data set, where $\BB{C}$ has a very large condition number and slowly decaying dominant eigenvalues.  In this case,  low-rank approximation methods with a small approximate rank $r$ may not be able to capture sufficient curvature information. We can see that even in this case  curvature information helps to avoid the stagnation experienced by \texttt{FISTA} and \texttt{ProxSVRG}. 
}{
Figures~1 and~2 show the suboptimality in objective versus the number of epochs\deleted{(passes over the full data set)} for different algorithms solving the elastic net problem. We can see that our method systematically outperforms the others in all setting\added{s}. Notably, for the australian dataset, all the first-order methods make almost no progress in the first 100 epochs, while our method can find a high\added{-}accuracy solution within tens of epochs, demonstrating a great benefit of second-order information. On the other hand, for a well-conditioned data set such as real-sim, \texttt{ProxSVRG} \replaced{is}{performs quite}  comparable to our method and even outperforms \texttt{Katyusha1}. This\deleted{highly} agrees with the theoretical time complexities summarized in Table~\ref{tab:comparison}.}

\section{Conclusions}

We have proposed and analyzed a novel second-order method for solving the elastic net problem. {By carefully exploiting} the {problem structure, we demonstrated that it is possible} to deal with the non-smooth objective and {to} efficiently inject curvature information into the optimization process without (significantly) increasing the computational cost per iteration. The combination of second-order information, fast iterative solvers, and a well-designed warm-start procedure results in a significant improvement of the total runtime complexity over popular first-order methods. An interesting direction for future research would be to go beyond the quadratic loss. We believe that the techniques developed in this work can be extended to more general settings, especially when the smooth part of the objective function is self-concordant.

\section*{Acknowledgments}
This research was sponsored in part by the Knut and Alice Wallenberg
Foundation and the Swedish Research Council.

\bibliographystyle{plain}
\bibliography{refs}
\newpage

\appendix

\section{Proof of Theorem~\ref{thrm:bounding:cond}}\label{apdx:thrm:bounding:cond}
First, by the definition of $L_i$, it holds that
\begin{align*}
	L_i
	=
		\Vert\He^{-1}\big(\av_i\av_i^\top+\reg_2\IM{}\big)\Vert_2
	\leq
		\av_i^\top \He^{-1}\av_i
		+
		\reg_2\Vert\He^{-1}\Vert_2
	\leq
		\av_i^\top \He^{-1}\av_i
		+ 
		1,
\end{align*}
where the last inequality follows since
\begin{align*}
	\reg_2\Vert\He^{-1}\Vert_2
	=
		\frac{\reg_2}{\Vert\He\Vert_2}
	=
		\frac{\reg_2}{\eig_1\left(\He\right)+\reg_2}
	\leq 1.
\end{align*}

The term $\av_i^\top \He^{-1}\av_i$ can be bounded as
\begin{align*}
	\av_i^\top\He^{-1}\av_i
	&=
		\av_i^\top\left(\He^{-1}\left(\BB{C}+\reg_2\IM{}\right)-\IM{}\right)\left(\BB{C}+\reg_2\IM{}\right)^{-1}\av_i
		+
		\av_i^\top\left(\BB{C}+\reg_2\IM{}\right)^{-1}\av_i
	\nonumber\\
	&\leq
		\left(\norm{\eig_1\left(\He^{-1}\left(\BB{C}+\reg_2\IM{}\right)-\IM{}\right)}+1\right)
		\av_i^\top\left(\BB{C}+\reg_2\IM{}\right)^{-1}\av_i.
\end{align*}
It can be verified that
\begin{align*}
	\norm{\eig_1\left(\He^{-1}\left(\BB{C}+\reg_2\IM{}\right)-\IM{}\right)}
	=
	\max\left\{
		\norm{\eig_1\left(\He^{-1}\left(\BB{C}+\reg_2\IM{}\right)\right)-1},
		\norm{\eig_d\left(\He^{-1}\left(\BB{C}+\reg_2\IM{}\right)\right)-1}
	\right\}	
	\leq 16,
\end{align*}
where we have used Lemma~\ref{lem:eig:bound} in the last step.
It follows that
\begin{align*}
	L_i
	\leq
		17\av_i^\top\left(\BB{C}+\reg_2\IM{}\right)^{-1}\av_i
		+
		1.
\end{align*}
Note that
\begin{align*}	
	\frac{1}{n}\sum_{i=1}^{n}\av_i^\top\left(\BB{C}+\reg_2\IM{}\right)^{-1}\av_i
	&=	
		\frac{1}{n}\sum_{i=1}^{n}\tr\left(\left(\BB{C}+\reg_2\IM{}\right)^{-1}\av_i\av_i^\top\right)
	\\
	&=	
		\tr\left(\left(\BB{C}+\reg_2\IM{}\right)^{-1}\BB{C}\right)
	\\
	&=
		\sum_{i=1}^{d}
		\frac{\eig_i}{\eig_i+\reg_2}
	= 
		\sdim_{\reg_2},
\end{align*}
which implies that
\begin{align*}
	\Lavg
	=
		\frac{1}{n}\sum_{i=1}^{n}L_i
	\leq
		17\sdim_{\reg_2}+1
	\leq 18\sdim_{\reg_2}.
\end{align*}
On the other hand, it was shown in~\cite[Theorem~4]{GOS16} that 
\begin{align*}
	\Lavg
	\leq
		\tr\left(\He^{-1/2}\left(\BB{C}+\reg_2\IM{}\right)\He^{-1/2}\right)
	=
	O\left(r\eig_r+\sum_{i>r}\eig_i+\reg_2 d\right).
\end{align*}
Finally, recall that $\mu=\eig_d\left(\He^{-1/2}\left(\BB{C}+\reg_2\IM{}\right)\He^{-1/2}\right)$. By Lemma~\ref{lem:eig:bound}, we readily have $\mu\geq c_1\reg_2$ for some constant $c_1$, which concludes the proof. 

\section{Some Useful Auxiliary Results}
To {facilitate} the analysis, we collect some useful inequalities regarding the Mahalanobis norm that \replaced{are used}{will be used extensively} in the subsequent proofs.
\begin{itemize}
\item Cauchy' inequality: $\InP{\xv}{\yv}\leq\Mnorm{\xv}{\He} \Mnorm{\yv}{\He^{-1}}$  for all $\xv, \yv \in \R^d$
\item Young's inequality: $\InP{\xv}{\yv}\leq \frac{\alpha\Mnorm{\xv}{\He}^2}{2}+\frac{\Mnorm{\yv}{\He^{-1}}^2}{2\alpha}$ for all $\xv, \yv \in \R^d$
\item Strongly convex inequality: $\frac{\mu}{2}\Mnorm{\xv-\xv^\star}{\He}^2\leq F\left
(\xv\right)-F\left(\xv^\star\right)$, where $\xv^\star=\argmin_{\xv}F\left(\xv\right)$.
\end{itemize}

\section{Proof of Lemma~\ref{keylem}}\label{apdx:keylem}
Before proving the lemma, we rewrite the sequences $\{\xv_k,\yv_k,\zv_k\}_{k\geq0}$ generated by Algorithm~\ref{alg:1} in the following recurrence form:
\begin{subequations}
\begin{align}
	\xv_{k+1} \label{eq:key:a}
	&= 
		\yv_k - \eta \, \gv_{k+1} 
	\\
	\tau\left(\yv_k-\zv_k\right) \label{eq:key:b}
	&=
		\xv_k-\yv_k
	\\
	\zv_{k+1} - \zv_k \label{eq:key:c}
	&=
		 \tau\left(\yv_k-\zv_k\right) - \frac{\tau}{\mu }\gv_{k+1}.
\end{align}
\end{subequations}
We also recall the following well-known three-point identity:
\begin{align}\label{eq:2p}
	\InP{\BB{a}-\BB{b}}{\BB{a}-\BB{c}}_{\He}
	=
		\frac{1}{2}		
		\Mnorm{\BB{a}-\BB{b}}{\He}^2
		+
		\frac{1}{2}		
		\Mnorm{\BB{a}-\BB{c}}{\He}^2
		-
		\frac{1}{2}		
		\Mnorm{\BB{b}-\BB{c}}{\He}^2,
\end{align}
which holds for any symmetric matrix $\He \succ \B{0}$ and any vectors $\BB{a}, \BB{b}, \BB{c} \in \R^d$.\\

We are now ready to prove Lemma~\ref{keylem}. From the definition of the Lyapunov function in~\eqref{eq:Lyp:def}, we have that
\begin{align}\label{eq:proof:keylem:1}
	\Lyp_{k+1}-\Lyp_k 
	&= 
		\FF{\xv_{k+1}}-\FF{\xv_k} 
		+
		\frac{\mu}{2}\Mnorm{\zv_{k+1}-\xv^\star}{\He}^2
		-
		\frac{\mu}{2}\Mnorm{\zv_{k}-\xv^\star}{\He}^2
	\nonumber\\
	&\mathop = \limits^{\mathrm{(a)}}
		\FF{\xv_{k+1}}-\FF{\xv_k}
		-
		\mu
		\InP{\zv_{k+1}-\zv_k}{\xv^\star - \zv_{k+1}}_{\He}
		-
		\frac{\mu}{2}\Mnorm{\zv_{k+1}-\zv_{k}}{\He}^2
	\nonumber\\
	&=
		\FF{\xv_{k+1}}-\FF{\xv_k}
		-
		\mu
		\InP{\zv_{k+1}-\zv_k}{\xv^\star - \zv_k}_{\He}
		+
		\frac{\mu}{2}\Mnorm{\zv_{k+1}-\zv_{k}}{\He}^2
	\nonumber\\
	&\mathop = \limits^{\mathrm{(b)}}
		\FF{\xv_{k+1}}-\FF{\xv_k}
		+
		\frac{\mu}{2}\Mnorm{\zv_{k+1}-\zv_{k}}{\He}^2
		+
		\InP{\tau\gv_{k+1}-\mu\tau\left(\yv_k-\zv_k\right)}{\xv^\star - \zv_k}_{\He}
	\nonumber\\
	&=
		\FF{\xv_{k+1}}-\FF{\xv_k}
		+
		\frac{\mu}{2}\Mnorm{\zv_{k+1}-\zv_{k}}{\He}^2
		+
		\tau
		\InP{\gv_{k+1}}{\xv^\star - \yv_k}_{\He}
	\nonumber\\
	&\hspace{3.25cm}	
		+
		\tau
		\InP{\gv_{k+1}}{\yv_k - \zv_k}_{\He}
		-\mu\tau	
		\InP{\yv_k-\zv_k}{\xv^\star - \zv_k}_{\He},
\end{align}
where the equality $\mathrm{(a)}$ follows from~\eqref{eq:2p}, and $\mathrm{(b)}$ follows from\eqref{eq:key:c}.
Using the identity~\eqref{eq:2p} again for the last term in~\eqref{eq:proof:keylem:1}, we obtain
\begin{align}\label{eq:proof:keylem:2}
	\Lyp_{k+1}-\Lyp_k 
	&= 
		\FF{\xv_{k+1}}-\FF{\xv_k}
		+
		\frac{\mu}{2}\Mnorm{\zv_{k+1}-\zv_{k}}{\He}^2
	\nonumber\\
	&\hspace{0.5cm}
		+
		\tau
		\InP{\gv_{k+1}}{\xv^\star - \yv_k}_{\He}	
		+
		\tau
		\InP{\gv_{k+1}}{\yv_k - \zv_k}_{\He}
		\nonumber\\
	&\hspace{0.5cm}
		-
		\frac{\mu\tau}{2}
		\Mnorm{\yv_k-\zv_k}{\He}^2
		-
		\frac{\mu\tau}{2}
		\Mnorm{\zv_k-\xv^\star}{\He}^2
		+
		\frac{\mu\tau}{2}
		\Mnorm{\yv_k-\xv^\star}{\He}^2
	\nonumber\\
	&= 
		\FF{\xv_{k+1}}-\FF{\xv_k}
		+
		\frac{\mu}{2}\Mnorm{\zv_{k+1}-\zv_{k}}{\He}^2
	\nonumber\\
	&\hspace{0.5cm}
		+
		\tau
		\InP{\gv_{k+1}}{\xv^\star - \yv_k}_{\He}	
		+
		\InP{\gv_{k+1}}{\xv_k - \yv_k}_{\He}
	\nonumber\\
	&\hspace{0.5cm}
		-
		\frac{\mu}{2\tau}
		\Mnorm{\xv_k-\yv_k}{\He}^2
		-
		\frac{\mu\tau}{2}
		\Mnorm{\zv_k-\xv^\star}{\He}^2
		+
		\frac{\mu\tau}{2}
		\Mnorm{\yv_k-\xv^\star}{\He}^2,
\end{align}
where the last equality follows from~\eqref{eq:key:b}. By noting that
\begin{align}\label{eq:proof:keylem:3}
	\frac{\mu}{2}\Mnorm{\zv_{k+1}-\zv_{k}}{\He}^2
	&=
		\frac{\mu \tau^2}{2}
		\Mnorm{\yv_k-\zv_k - \frac{1}{\mu}\gv_{k+1}}{\He}^2
	\nonumber\\
	&=
		\frac{\mu \tau^2}{2}\left(
			\Mnorm{\yv_k-\zv_k}{\He}^2
			-
			\frac{2}{\mu}
			\InP{\gv_{k+1}}{\yv_k-\zv_k}_{\He}
			+
			\frac{1}{\mu^2}
			\Mnorm{\gv_{k+1}}{\He}^2
		\right)
	\nonumber\\
	&=
		\frac{\mu }{2}
			\Mnorm{\xv_k-\yv_k}{\He}^2
			-
			\tau
			\InP{\gv_{k+1}}{\xv_k-\yv_k}_{\He}
			+
			\frac{\tau^2}{2\mu}
			\Mnorm{\gv_{k+1}}{\He}^2,
\end{align}
where we have used~\eqref{eq:key:b} in the last equality. Then, combining~\eqref{eq:proof:keylem:2} and~\eqref{eq:proof:keylem:3} yields
\begin{align}\label{eq:proof:keylem:4}
	\Lyp_{k+1}-\Lyp_k 
	&=
		\FF{\xv_{k+1}}-\FF{\xv_k}
		+
		\frac{\tau^2}{2\mu}
		\Mnorm{\gv_{k+1}}{\He}^2
	\nonumber\\
	&\hspace{0.4cm}		
		+
		\tau
		\InP{\gv_{k+1}}{\xv^\star - \yv_k}_{\He}	
		+
		\left(1-\tau\right)
			\InP{\gv_{k+1}}{\xv_k-\yv_k}_{\He}	
	\nonumber\\
	&\hspace{0.4cm}
		+
		\left(\frac{\mu}{2}-\frac{\mu}{2\tau}\right)
		\Mnorm{\xv_k-\yv_k}{\He}^2
		-
		\frac{\mu\tau}{2}
		\Mnorm{\zv_k-\xv^\star}{\He}^2
		+
		\frac{\mu\tau}{2}
		\Mnorm{\yv_k-\xv^\star}{\He}^2.
\end{align}

We now pay attention to the term $\FF{\xv_{k+1}}$. Let $\uv=\left(1-\tau\right)\xv_k+\tau\xv^\star$, then it follows from the strong convexity of $F$ that
\begin{align}\label{eq:proof:keylem:5}
	\FF{\uv} 
	&\leq 
		\left(1-\tau\right) \FF{\xv_k} 
		+
		\tau \FF{\xv^\star}
		-
		\frac{\tau\left(1-\tau\right)\mu}{2}\Mnorm{\xv_k-\xv_\star}{\He}^2.	
\end{align}
In addition, we have that 
\begin{align}\label{eq:proof:keylem:6}
	&\frac{1}{2\eta}\left(\Mnorm{\xv_{k+1}-\uv}{\He}^2-\Mnorm{\yv_k-\uv}{\He}^2\right)
	\nonumber\\	
	&\hspace{0.5cm}
	=
		\frac{1}{2\eta}\left(
			\Mnorm{\xv_{k+1}-\yv_k}{\He}^2
			+
			2\InP{\xv_{k+1}-\yv_k}{\yv_k-\uv}_{\He}
		\right)
	\nonumber\\
	&\hspace{0.5cm}
	=
		\frac{\eta}{2}\Mnorm{\gv_{k+1}}{\He}^2
		-
		\left(1-\tau\right)
		\InP{\gv_{k+1}}{\yv_k-\xv_k}_{\He}
		-
		\tau
		\InP{\gv_{k+1}}{\yv_k-\xv^\star}_{\He},
\end{align}
where we have used~\eqref{eq:key:a} and the fact that $\yv_k-\uv=\left(1-\tau\right)\left(\yv_k-\xv_k\right)+\tau \left(\yv_k - \xv^\star\right)$.

With these observations, we are now ready to bound $\FF{\xv_{k+1}}$. In particular, by invoking Lemma~\ref{lem:svrg:progress} with $\yv=\yv_t$, $\xv^+=\xv_{t+1}$, $\uv = \left(1-\tau\right)\xv_k+\tau\xv^\star$, $\vv=\vv_t$, $\BB{\xi}=\BB{\xi}_k$, $\epsilon=\epsilon_k$, together with~\eqref{eq:proof:keylem:5}--\eqref{eq:proof:keylem:6}, it follows that
\begin{align}\label{eq:proof:keylem:7}
	\FF{\xv_{k+1}}
	&\leq
		\left(1-\tau\right) \FF{\xv_k} 
		+
		\tau \FF{\xv^\star}
		-
		\frac{\tau\left(1-\tau\right)\mu}{2}\Mnorm{\xv_k-\xv_\star}{\He}^2
	\nonumber\\
	&\hspace{0.4cm}
		-
		\frac{\eta}{2}\Mnorm{\gv_{k+1}}{\He}^2
		+
		\left(1-\tau\right)
		\InP{\gv_{k+1}}{\yv_k-\xv_k}_{\He}
		+
		\tau
		\InP{\gv_{k+1}}{\yv_k-\xv^\star}_{\He}
	\nonumber\\
	&\hspace{0.4cm}
		-\frac{\mu}{2}
		\Mnorm{\yv_k-\left(1-\tau\right)\xv_k-\tau\xv^\star}{\He}^2
	+
	\InP{\Nvar_k}{\left(1-\tau\right)\xv_k+\tau\xv^\star-\xv_{k+1}}
	\nonumber\\
	&\hspace{0.4cm}
	+
	\frac{1}{\eta}
	\InP{\BB{\xi}_k}{\left(1-\tau\right)\xv_k+\tau\xv^\star-\xv_{k+1}}
	+\epsilon_k.
\end{align}
Substituting~\eqref{eq:proof:keylem:7} into~\eqref{eq:proof:keylem:2} and rearranging the terms to obtain
\begin{align}\label{eq:proof:keylem:8}
	\Lyp_{k+1}-\Lyp_k 
	&\leq	
		-\tau\overbrace{\big(
		\FF{\xv_k}-\FF{\xv^\star}
		+
		\frac{\mu}{2}
		\Mnorm{\zv_k-\xv^\star}{\He}^2
		\big)}^{\Lyp_k}
		+
		\left(\frac{\tau^2}{2\mu}-\frac{\eta}{2}\right)
		\Mnorm{\gv_{k+1}}{\He}^2	
	\nonumber\\
	&\hspace{0.4cm}	
		+
		\left(\frac{\mu}{2}-\frac{\mu}{2\tau}\right)
		\Mnorm{\xv_k-\yv_k}{\He}^2
		+
		\frac{\mu\tau}{2}
		\Mnorm{\yv_k-\xv^\star}{\He}^2
	\nonumber\\
	&\hspace{0.4cm}
		+
		\underbrace{\InP{\Nvar_k}{\left(1-\tau\right)\xv_k+\tau\xv^\star-\xv_{k+1}}}_{T_3}
	\nonumber\\
	&\hspace{0.4cm}
		-\frac{\mu}{2}\underbrace{\big(
		\tau\left(1-\tau\right)\Mnorm{\xv_k-\xv_\star}{\He}^2
		+
		\Mnorm{\yv_k-\left(1-\tau\right)\xv_k-\tau\xv^\star}{\He}^2
		\big)}_{T_2}
	\nonumber\\
	&\hspace{0.4cm}
	+
	\frac{1}{\eta}
	\underbrace{
	\InP{\BB{\xi}_k}{\left(1-\tau\right)\xv_k+\tau\xv^\star-\xv_{k+1}}
	}_{T_1}
	+\epsilon_k.
\end{align}

We next bound the term $T_1$. By adding and subtracting the term $\yv_k$ yields
\begin{align*}
	T_1
	&=
		\InP{\BB{\xi}_k}{\xv_k-\yv_k + \tau\left(\xv^\star-\xv_k\right) + \yv_k-\xv_{k+1}}
	\nonumber\\
	&=
		\InP{\BB{\xi}_k}{\tau\left(\yv_k-\zv_k\right) + \tau\left(\xv^\star-\xv_k\right) + \eta\gv_{k+1}}
		\nonumber\\
	&=		
		\tau	
		\InP{\BB{\xi}_k}{\xv^\star-\zv_k}
		+
		\tau
		\InP{\BB{\xi}_k}{\yv_k-\xv_k}
		+
		\eta
		\InP{\BB{\xi}_k}{\gv_{k+1}},
\end{align*}
where we have used~\eqref{eq:key:a} and~\eqref{eq:key:b}.
By Young's inequality, we have
\begin{align}\label{eq:T1:1}
	\tau
	\InP{\BB{\xi}_k}{\yv_k-\xv_k}
	&\leq
		\frac{\Mnorm{\BB{\xi}_k}{\He^{-1}}^2}{2}
		+
		\frac{\tau^2\Mnorm{\xv_k-\yv_k}{\He}^2}{2}
	\\\label{eq:T1:2}
	\InP{\BB{\xi}_k}{\gv_{k+1}}
	&\leq
		\frac{\Mnorm{\BB{\xi}_k}{\He^{-1}}^2}{\eta}
		+
		\frac{\eta\Mnorm{\gv_{k+1}}{\He}^2}{4}.
\end{align}
For $T_2$, we have
\begin{align}\label{eq:proof:keylem:10}
	T_2
	&=	
		\tau\left(1-\tau\right)\Mnorm{\xv_k-\xv_\star}{\He}^2
		+
		\Mnorm{\yv_k-\xv_k + \tau\left(\xv_k-\xv^\star\right)}{\He}^2	
	\nonumber\\
	&=
		\Mnorm{\yv_k-\xv_k}{\He}^2	
		+
		2\tau\InP{\yv_k-\xv_k }{\xv_k-\xv^\star}_{\He}
		+
		\tau \Mnorm{\xv_k-\xv^\star}{\He}^2
	\nonumber\\
	&=
		 \left(1-\tau\right)\Mnorm{\yv_k-\xv_k}{\He}^2 + \tau \Mnorm{\yv_k-\xv^\star}{\He}^2.
\end{align}
Thus, combining~\eqref{eq:proof:keylem:8}--\eqref{eq:proof:keylem:10} and using the fact that $\Mnorm{\BB{\xi}_k}{\He^{-1}} \leq \sqrt{2\eta\epsilon_k}$ yield
\begin{align}\label{eq:proof:keylem:11}
	\Lyp_{k+1} 
	&\leq	
		\left(1-\tau\right)\Lyp_k
		+
		\left(\frac{\tau^2}{2\mu}-\frac{\eta}{4}\right)
		\Mnorm{\gv_{k+1}}{\He}^2	
		+
		\left(\frac{\tau^2}{2\eta}+\frac{\tau\mu}{2}-\frac{\mu}{2\tau}\right)
		\Mnorm{\xv_k-\yv_k}{\He}^2
	\nonumber\\
	&\hspace{1.cm}
		+
		T_3
		+
		\frac{\tau}{\eta}	
		\InP{\BB{\xi}_k}{\xv^\star-\zv_k}
		+
		4\epsilon_k.
\end{align}
Invoking Lemma~\ref{lem:var:grad} with $\uv_k=\left(1-\tau\right)\xv_k+\tau\xv^\star$ and then applying Lemma~\ref{lem:bouding:var}, we obtain
\begin{align}\label{eq:proof:keylem:11b}
	\E T_3
	\leq
		\frac{3\eta \Lavg^2}{b}	
		\E\Mnorm{\xv_k-\yv_k}{\He}^2
		+
		\frac{12\eta \Lavg}{b}	\EX{
			\FF{\xv_k} - \FF{\xv^\star}
			+
			\FF{\xTd_s} - \FF{\xv^\star}
		}	
		+
		\epsilon_k.	
\end{align}
Therefore, by taking the expectation on both sides of~\eqref{eq:proof:keylem:11} and using~\eqref{eq:proof:keylem:11b}, we obtain
\begin{align}\label{eq:proof:keylem:12}
	\E\Lyp_{k+1} 
	&\leq	
		\left(1-\tau\right)\E\Lyp_k
		+
		\left(\frac{\tau^2}{2\mu}-\frac{\eta}{4}\right)
		\E\Mnorm{\gv_{k+1}}{\He}^2	
	\nonumber\\
	&\hspace{0.4cm}
		+
		\left(\frac{\tau^2}{2\eta}+\frac{\tau\mu}{2}-\frac{\mu}{2\tau} + \frac{3\eta \Lavg^2}{b}\right)
		\E\Mnorm{\xv_k-\yv_k}{\He}^2
	\nonumber\\
	&\hspace{0.4cm}
		+
		\frac{12\eta \Lavg}{b}
		\EX{
			\FF{\xv_k} - \FF{\xv^\star}
			+
			\FF{\xTd_s} - \FF{\xv^\star}
		}
		+\frac{\tau}{\eta}	
		\E\InP{\BB{\xi}_k}{\xv^\star-\zv_k}
		+
		5\epsilon_k.
\end{align}
By choosing $\eta = \frac{1}{\Lavg}$, $\tau=\sqrt{\frac{\mu}{2\Lavg}}$, and $b \geq 60\sqrt{\frac{\Lavg}{\mu}}$, it is readily verified that the second and third terms on the right-hand side of~\eqref{eq:proof:keylem:12} become nonpositive, which concludes the proof.

\section{Proof of Theorem~\ref{thrm:main}}\label{apdx:thrm:main}
To begin with, let $e_t=\EX{\FF{\xv_t} - \FF{\xv^\star}}$ and $\tilde{e}_s=\EX{\FF{\xTd_s} - \FF{\xv^\star}}$, then by applying inequality~\eqref{eq:keylem} recursively, we obtain
\begin{align*}
	\E\Lyp_{k} 
	&\leq	
		\left(1-\tau\right)^k\Lyp_0		
		+
		\frac{\tau}{5}
		\sum_{t=0}^{k-1}
		\left(1-\tau\right)^{k-t-1}		
		\left(e_t+\tilde{e}_s\right)		
	\nonumber\\
	&\hspace{1cm}
		+	
		\E		
		\sum_{t=0}^{k-1}
		\left(1-\tau\right)^{k-t-1}
		\left(
			\Lavg\tau	
			\Mnorm{\BB{\xi}_t}{\He^{-1}}
			\Mnorm{\zv_t-\xv^\star}{\He}
			+
			5\epsilon_t
		\right).
\end{align*}
It can be verified that  
\begin{align}\label{eq:geo:sum}
	\sum_{t=0}^{k-1}
	\left(1-\tau\right)^{k-t-1}
	=
		\sum_{t=0}^{k-1}
		\left(1-\tau\right)^{t}
	=
		\frac{1-\left(1-\tau\right)^{k}}{\tau}
	\leq
		\frac{1}{\tau}.
\end{align}
Thus, from~\eqref{eq:geo:sum} and the facts that $\tau=\sqrt{\frac{\mu}{2\Lavg}}$, $\Mnorm{\BB{\xi}_t}{\He^{-1}}\leq\sqrt{2\epsilon_t/\Lavg}$, it holds that
\begin{align}\label{eq:proof:main:thrm:main}
	\E\Lyp_{k} 
	&\leq
		\left(1-\tau\right)^k\Lyp_0		
		+	
		5
		\sum_{t=0}^{k-1}
		\left(1-\tau\right)^{k-t-1}				
			\epsilon_t
		+
		\frac{\tau}{5}
		\sum_{t=0}^{k-1}
		\left(1-\tau\right)^{k-t-1}		
		e_t		
		+
		\frac{\tilde{e}_s}{5}		
	\nonumber\\
	&\hspace{1.5cm}
		+
		\sum_{t=0}^{k-1}
		\left(1-\tau\right)^{k-t-1}
			\sqrt{\mu\epsilon_t}
			\E\Mnorm{\zv_t-\xv^\star}{\He}
	\nonumber\\
	&\leq
		\left(1-\tau\right)^k
		\left(\Lyp_0+B_k\right)
		+
		\sum_{t=0}^{k-1}
			\left(1-\tau\right)^{k-t-1}
			\sqrt{\mu\epsilon_t}
			\E\Mnorm{\zv_t-\xv^\star}{\He},
\end{align}
where 
\begin{align*}
	B_k
	=
		5\sum_{t=0}^{k-1}
		\left(1-\tau\right)^{-t-1}				
			\epsilon_t
		+
		\frac{\tau}{5}
		\sum_{t=0}^{k-1}
		\left(1-\tau\right)^{-t-1}		
		e_t	
		+		
		\left(1-\tau\right)^{-k}
		\frac{\tilde{e}_s}{5}.
\end{align*}
By the definition of $\Lyp_k$, we have $\frac{\mu}{2}\Mnorm{\zv_k-\xv^\star}{\He}^2\leq \Lyp_k$, hence, it follows from~\eqref{eq:proof:main:thrm:main} that
\begin{align*}
	\E\Mnorm{\zv_k-\xv^\star}{\He}^2
	&\leq
		\frac{2}{\mu}
		\left(1-\tau\right)^k
		\left(\Lyp_0+B_k\right)	
		+	
		2
		\sum_{t=0}^{k-1}
		\sqrt{\frac{\epsilon_t}{\mu}}
		\left(1-\tau\right)^{k-t-1}
		\E\Mnorm{\zv_t-\xv^\star}{\He}.	
\end{align*}
Multiplying both sides of the above inequality by $\left(1-\tau\right)^{-k}$ yields
\begin{align*}
	\left(1-\tau\right)^{-k}
	\E\Mnorm{\zv_k-\xv^\star}{\He}^2
	&\leq
		\frac{2}{\mu}
		\left(\Lyp_0+B_k\right)
		+	
		2
		\sum_{t=0}^{k-1}
		\sqrt{\frac{\epsilon_t}{\mu}}
		\left(1-\tau\right)^{-t/2-1}
		\left(1-\tau\right)^{-t/2}
		\E\Mnorm{\zv_t-\xv^\star}{\He}								
\end{align*}
Define $\alpha_t=2\sqrt{\frac{\epsilon_t}{\mu}}\left(1-\tau\right)^{-t/2-1}$ and $u_t=\left(1-\tau\right)^{-t/2}\E\Mnorm{\zv_t-\xv^\star}{\He}$, then we can write the previous inequality as
\begin{align*}
	u_k^2
	\leq
		\frac{2}{\mu}
		\left(\Lyp_0+B_k\right)
		+
		\sum_{t=0}^{k-1}		
		\alpha_t
		u_t
	\leq
		\frac{2}{\mu}
			\left(\Lyp_0+B_k\right)
			+
			\alpha_0
			u_0
			+
			\sum_{t=1}^{k}		
			\alpha_t
			u_t,
\end{align*}
where we have separated the term $\alpha_0 u_0$ from the sum and added the positive term  $\alpha_k u_k$ in the last step.
Note that $\alpha_0=\frac{2}{1-\tau}\sqrt{\frac{\epsilon_0}{\mu}}\leq\frac{2}{1-\tau}\sqrt{\frac{\Lyp_0}{\mu}}$ and $u_0=\Mnorm{\xv_0-\xv^\star}{\He}\leq \sqrt{\frac{2}{\mu}\Lyp_0}$, hence, $u_k^2$ can be further bounded by
\begin{align*}
	u_k^2
	&\leq
		\frac{2}{\mu}
		\left(c_1\Lyp_0+B_k\right)
		+
		\sum_{t=1}^{k}		
		\alpha_t
		u_t,
\end{align*}
where  $c_1=1+\frac{\sqrt{2}}{1-\tau}$.	It is readily verified that \added{$\{ B_k\}_k$ is an increasing sequence} and that   $\{S_k\triangleq\frac{2}{\mu}\left(c_1\Lyp_0+B_k\right)\}_k$ is an increasing sequence satisfying $S_0> u_0^2$. Therefore, by invoking Lemma~\ref{lem:iter:recurs} with $u_t$, $\alpha_k$, and $S_k$, we have for any $k\geq 1$ that
\begin{align*}
	\left(1-\tau\right)^{-k/2}
	\E\Mnorm{\zv_k-\xv^\star}{\He}
	&\leq
		\frac{1}{2}
		\sum_{t=1}^{k}
		\alpha_t
		+
		\left(
			\frac{2c_1}{\mu}
			\Lyp_0
			+
			\frac{2}{\mu}
			B_k
			+
			\left(\frac{1}{2}\sum_{t=1}^{k}\alpha_t\right)^2
		\right)^{1/2}
	\nonumber\\
	&\leq
		\sum_{t=1}^{k}
		\alpha_t
		+
		\sqrt{\frac{2c_1}{\mu}}
		\sqrt{\Lyp_0}
		+
		\sqrt{\frac{2}{\mu}}
		\sqrt{B_k},
\end{align*}
where we have used  $\sqrt{a+b}\leq \sqrt{a}+\sqrt{b}$ for any $a,b \geq 0$. Thus, for any $t\in\{0,1,\ldots,k\}$, it holds that 
\begin{align*}	
	\E\Mnorm{\zv_t-\xv^\star}{\He}
	&\leq
		\left(1-\tau\right)^{t/2}
		\left(
		\sum_{i=1}^{t}
		\alpha_i
		+
		\sqrt{\frac{2c_1}{\mu}}
		\sqrt{\Lyp_0}
		+
		\sqrt{\frac{2}{\mu}}
		\sqrt{B_t}
		\right)
	\\
	&\leq
		\left(1-\tau\right)^{t/2}
		\left(
		\sum_{i=1}^{k}
		\alpha_i
		+
		\sqrt{\frac{2c_1}{\mu}}
		\sqrt{\Lyp_0}
		+
		\sqrt{\frac{2}{\mu}}
		\sqrt{B_k}
		\right).
\end{align*}
Having upper bounds of $\E\Mnorm{\zv_t-\xv^\star}{\He}$, we can now substitute them into~\eqref{eq:proof:main:thrm:main} to get
\begin{align*}
	\E\Lyp_k
	&\leq
	\left(1-\tau\right)^k
	\left(
		\Lyp_0	
		+
		B_k
	\right)
	+
	\sum_{t=0}^{k-1}
	\left(1-\tau\right)^{k-t/2-1}
		\sqrt{\mu\epsilon_t}
		\left(
			\sum_{t=1}^{k}
			\alpha_t
			+
			\sqrt{\frac{2c_1}{\mu}}
			\sqrt{\Lyp_0}
			+
			\sqrt{\frac{2}{\mu}}
			\sqrt{B_k}
			\right)
	\nonumber\\
	&=
	\left(1-\tau\right)^k
		\left(
			\Lyp_0	
			+
			B_k	
		+
		\frac{\mu}{2}
		\sum_{t=0}^{k-1}
			\alpha_t
			\left(
				\sum_{t=1}^{k}
				\alpha_t
				+
				\sqrt{\frac{2c_1}{\mu}}
				\sqrt{\Lyp_0}
				+
				\sqrt{\frac{2}{\mu}}
				\sqrt{B_k}
				\right)	
		\right)
	\nonumber\\
	&\leq
	\left(1-\tau\right)^k
		\left(
			\Lyp_0	
			+
			B_k	
			+
			\frac{\mu}{2}
			\sum_{t=0}^{k}
				\alpha_t
				\left(
					\sum_{t=0}^{k}
					\alpha_t
					-
					\alpha_0
					+
					\sqrt{\frac{2c_1}{\mu}}
					\sqrt{\Lyp_0}
					+
					\sqrt{\frac{2}{\mu}}
					\sqrt{B_k}
				\right)	
		\right),
\end{align*}
where the equality follows from the definition of $\alpha_t$, and in the last step, we have added a positive term $\alpha_k$ in the first sum as well as added and subtracted $\alpha_0$ in the second sum.
Note that the function $\sqrt{1+x}-x$ is decreasing on the interval~$\left[0,\infty\right)$, it follows that
\begin{align*}
	\sqrt{\frac{2c_1}{\mu}}\sqrt{\Lyp_0}-\alpha_0
		=
			\left(\sqrt{1+\frac{\sqrt{2}}{1-\tau}}-\frac{\sqrt{2}}{1-\tau}\right)
			\sqrt{\frac{2}{\mu}}\sqrt{\Lyp_0}
		\leq
			\sqrt{\frac{2}{\mu}}\sqrt{\Lyp_0}.
\end{align*}
We thus have
\begin{align*}
	\E\Lyp_k
	&\leq
		\left(1-\tau\right)^k
		\left(
			\Lyp_0	
			+
			B_k	
			+
			\frac{\mu}{2}
			\sum_{t=0}^{k}
				\alpha_t
				\left(
					\sum_{t=0}^{k}
					\alpha_t
					+
					\sqrt{\frac{2}{\mu}}
					\sqrt{\Lyp_0}
					+
					\sqrt{\frac{2}{\mu}}
					\sqrt{B_k}
				\right)	
		\right)
	\nonumber\\
	&\leq
	\left(1-\tau\right)^k
		\left(
			\sqrt{\Lyp_0}
			+						
		\sqrt{\frac{\mu}{2}}
		\sum_{t=0}^{k}
			\alpha_t
		+
		\sqrt{B_k}			
		\right)^2
	\nonumber\\
	& \leq
	\left(1-\tau\right)^k\left(
		20\Lyp_0
		+						
		10\mu
		\left(\sum_{t=0}^{k}\alpha_t\right)^2
		+
		\frac{10}{9}B_k
		\right),			
\end{align*}
 where in the last step, we have applied the inequality $\left(a+b\right)^2\leq \left(1+\beta\right)a^2+\left(1+\frac{1}{\beta}\right)b^2$ twice with $\beta=9$ and $\beta=1$, respectively.
Using the definitions of $B_k$ and $\alpha_t$, the above inequality can be rewritten as
\begin{align}\label{eq:proof:main:thrm:2}
	\E\Lyp_k
		&\leq
		20\left(1-\tau\right)^k
		\Lyp_0
		+
		\frac{50}{9}
		\sum_{t=0}^{k-1}
		\left(1-\tau\right)^{k-t-1}				
			\epsilon_t
		+
		\frac{2\tau}{9}
		\sum_{t=0}^{k-1}
		\left(1-\tau\right)^{k-t-1}		
		e_t	
		+
		\frac{2}{9}
		\tilde{e}_s
	\nonumber\\
	&\hspace{0.5cm}
	+						
	40
	\left(1-\tau\right)^k
	\left(\sum_{t=0}^{k}\sqrt{\epsilon_t}\left(1-\tau\right)^{-t/2-1}\right)^2.				
\end{align}
With our choice of $\epsilon_t$, it follows that
\begin{align*}
		\sum_{t=0}^{k}\sqrt{\epsilon_t}\left(1-\tau\right)^{-t/2-1}
	\leq
		\frac{1}{1-\tau}
		\sum_{t=0}^{k}
			\left(\sqrt{\frac{1-\rho}{1-\tau}}\,\right)^t
			\sqrt{\Lyp_0}
%
	\leq
		\frac{\left(\sqrt{\frac{1-\rho}{1-\tau}}\,\right)^{k+1}}{\sqrt{1-\rho}-\sqrt{1-\tau}}
		\sqrt{\frac{\Lyp_0}{1-\tau}}.
\end{align*}
Thus, the last term on the right-hand side of~\eqref{eq:proof:main:thrm:2} can be bounded by
\begin{align*}
	\frac{40\left(1-\rho\right)}{\left(\sqrt{1-\rho}-\sqrt{1-\tau}\right)^2\left(1-\tau\right)^2}	
	\left(1-\rho\right)^{k}\Lyp_0
	\leq
	\frac{160\left(1-\rho\right)}{\left(\tau-\rho\right)^2\left(1-\tau\right)^2}
	\left(1-\rho\right)^{k}\Lyp_0,
\end{align*}
where the last inequality follows since the function $\sqrt{1-x}+x/2$ is decreasing on the interval~$\left[0,1\right]$. 
Similarly, for the second term in~\eqref{eq:proof:main:thrm:2}, we have 
\begin{align*}
	\sum_{t=0}^{k-1}
		\left(1-\tau\right)^{k-t-1}				
		\epsilon_t
	\leq
	\left(1-\tau\right)^{k-1}
	\sum_{t=0}^{k-1}
	\left(\frac{1-\rho}{1-\tau}\right)^t
	\Lyp_0
	\leq
	\frac{\left(1-\rho\right)^k\Lyp_0}{\tau-\rho}.	
\end{align*}
Thus, we can further bound $\E\Lyp_k$ as
\begin{align*}
	\E\Lyp_k
		&\leq
		c_2\left(1-\rho\right)^k
		\Lyp_0
		+						
		\frac{2\tau}{9}
		\sum_{t=0}^{k-1}
		\left(1-\rho\right)^{k-t-1}		
		e_t	
		+
		\frac{2}{9}
		\tilde{e}_s,				
\end{align*}
where $c_2=20+\frac{160\left(1-\rho\right)}{\left(\tau-\rho\right)^2\left(1-\tau\right)^2}+\frac{50}{9\left(\tau-\rho \right)}$ and we have also used the fact that $1-\tau < 1-\rho$.

Now, if we let $\delta_k$ be the right-hand side of the above inequality, then it holds that
\begin{align*}
	\delta_k 
	&=
		c_2\left(1-\rho\right)^k
		\Lyp_0
		+						
	\frac{2\tau}{9}
		\sum_{t=0}^{k-1}
		\left(1-\rho\right)^{k-t-1}		
		e_t	
		+
		\frac{2}{9}
		\tilde{e}_s
	\nonumber\\
	&=
	\left(1-\rho\right)
	\left(
			c_2\left(1-\rho\right)^{k-1}
			\Lyp_0
			+						
			\frac{2\tau}{9}
			\sum_{t=0}^{k-2}
			\left(1-\rho\right)^{k-t-2}		
			e_t	
			+
			\frac{2}{9}
			\tilde{e}_s
		\right)	
		+
		\frac{2\tau}{9}		
		e_{k-1}	
		+
		\frac{2\rho}{9}
		\tilde{e}_s
	\nonumber\\
	&= 
		\left(1-\rho\right)\delta_{k-1}
		+
		\frac{2\tau}{9}	
		\delta_{k-1}
		+
		\frac{2\rho}{9}
		\Lyp_0		
	\nonumber\\
	&=
		\left(1-\frac{7\rho}{9}\right)\delta_{k-1}
		+
		\frac{2\rho}{9}
		\Lyp_0
	\nonumber\\
	&\leq
		\left(1-\frac{3\rho}{4}\right)\delta_{k-1}
		+
		\frac{2\rho}{9}
		\Lyp_0,	
\end{align*}
where the third equality follows since by definition, $e_{k-1 }\leq \E\Lyp_{k-1} \leq\delta_{k-1}$, and since $\tilde{e}_s=e_0\leq \Lyp_0$, and the last inequality follows since $\rho=0.9\tau$.
By applying the above inequality recursively, we obtain
\begin{align*}
	\E\Lyp_{T}\leq\delta_T 
	&\leq
		\left(1-\frac{3\rho}{4}\right)^{T}\delta_{0}
		+
		\sum_{t=0}^{T-1}
		\left(1-\frac{3\rho}{4}\right)^{T-t-1}
		\frac{2\rho}{9}
		\Lyp_0,		
	\nonumber\\
	&\leq
		\left(
			c_3\left(1-\frac{3\rho}{4}\right)^{T}	
			+
			\frac{8}{27}
		\right)
		\Lyp_0,
\end{align*}
where $c_3=c_2+2/9$. Therefore, when $T\geq\frac{4}{3\rho}\log \left(27c_3\right)$, it holds that 
\begin{align*}
	\E\Lyp_{T}\leq\delta_T 
	\leq
		\left(\frac{1}{27}+\frac{8}{27}\right)\Lyp_0
	=	
		\frac{1}{3}\Lyp_0.
\end{align*}
Finally, using the definition of $\Lyp_T$, we obtain
\begin{align*}
	\EX{\FF{\tilde{\BB{x}}_{s+1}}-\FF{\xv^\star}}
	\leq
		\E\Lyp_T
	\leq
		\frac{1}{3}\Lyp_0
	\leq
		\frac{2}{3}
		\left(\FF{\tilde{\BB{x}}_{s}}-\FF{\xv^\star}\right),
\end{align*}
completing the proof.

\section{Proof of Proposition~\ref{prop:ws}}\label{apdx:prop:ws}
Recall that 
\begin{align}\label{eq:def:z0}
	\zv_0=\Prox{\gamma h}{}{\xv_k-\frac{\gamma}{\eta}\He\left(\xv_k-\uv_{k-1}\right)},
\end{align}
which can be seen as one step of the proximal gradient method applied to $	p\left(\zv,\uv_{k-1} \right)$ starting at the current $\xv_k$ with
\begin{align*}
	\uv_k&=\yv_k - \eta\He^{-1} \vv_k
	\\
	p\left(\zv,\uv \right)&\triangleq\hF{\zv}+\frac{1}{2\eta} \Mnorm{\zv-\uv}{\He}^2.
\end{align*}
By \added{the} optimality condition of $\zv_0$, we have
\begin{align*}
	\frac{1}{\gamma}
	\left[
		-\zv_0
		+
		\xv_k
		-\frac{\gamma}{\eta}
		\He\left(\xv_k-\uv_{k-1} \right)
	\right]
	\in
	\partial\hF{\zv_0}.
\end{align*}
Since $\partial p\left(\zv_0,\uv_{k-1}\right)=\partial\hF{\zv_0}+\frac{1}{\eta}\He\left(\zv_0-\uv_{k-1} \right)$, it follows that
\begin{align*}
	\BB{\zeta}_k
	\triangleq
	\frac{1}{\gamma}
	\left(\xv_k-\zv_0\right)
	+	
	\frac{1}{\eta}\He\left(\zv_0-\xv_k\right)
	\in
	\partial p\left(\zv_0,\uv_{k-1}\right).
\end{align*}
We see that $\BB{\zeta}_k$ is independent of the second argument of $p\left(\zv_0,\cdot\right)$, hence we also have $\BB{\zeta}_k \in \partial p\left(\zv_0,\uv_{k}\right)$.
Since $p\left(\zv,\cdot\right)$ is  strongly convex in $\zv$ with parameter
 $\eig_r\left(\He\right)/\eta$ and $\BB{\zeta}_k \in \partial p\left(\zv_0,\uv_{k}\right)$, it holds that
\begin{align}\label{eq:proof:prop:ws:1}
	p\left(\xv_{k+1}^\star,\uv_k \right)
	&\geq
		p\left(\zv_0,\uv_k \right)
		+
		\InP{\BB{\zeta}_k}{\xv_{k+1}^\star-\zv_0}
		+
		\frac{\eig_r\left(\He\right)}{2\eta}
		\Vnorm{\xv_{k+1}^\star-\zv_0}_2^2
	\nonumber\\
	&\geq 
		p\left(\zv_0,\uv_k \right)
		-
		\frac{\eta}{2\eig_r\left(\He\right)}
		\Vnorm{\BB{\zeta}_k}^2,
\end{align}
where the last inequality follows from Young's inequality. We next bound $\Vnorm{\BB{\zeta}_k}^2$ \replaced{via}{, we have}
\begin{align}\label{eq:proof:prop:ws:2}
	\Vnorm{\BB{\zeta}_k}^2
	&=
		\Vnorm{
			\frac{1}{\gamma}
			\left(\xv_k-\zv_0\right)
			+	
			\frac{1}{\eta}\He\left(\zv_0-\xv_k\right)
		}_2^2
	\nonumber\\
	&=	
		\frac{1}{\gamma^2}
		\Vnorm{		
			\xv_k-\zv_0
		}_2^2						
		-
		\frac{2}{\gamma\eta}
		\InP{\xv_k-\zv_0}{\He\left(\xv_k-\zv_0\right)}
		+
		\frac{1}{\eta^2}
		\Vnorm{
			\He\left(\zv_0-\xv_k\right)
		}_2^2						
	\nonumber\\
	&\leq 
		\frac{1}{\gamma^2}	
		\Vnorm{
			\xv_k-\zv_0
		}_2^2,								
\end{align}
where the last step follows from the definition of $\gamma$ and the fact that 
\begin{align*}
	\InP{\xv_k-\zv_0}{\He\left(\xv_k-\zv_0\right)}
	=		
		\left(\xv_k-\zv_0\right)^\top\He\He^{-1}\He\left(\xv_k-\zv_0\right)
	\geq 			
	\frac{1}{\eig_1\left(\He\right)}
	\Vnorm{\He\left(\xv_k-\zv_0\right)}_2^2.								
\end{align*}
Thus, combining~\eqref{eq:proof:prop:ws:1} and~\eqref{eq:proof:prop:ws:2} yields
\begin{align*}
	p\left(\zv_0,\uv_k \right)-	p\left(\xv_{k+1}^\star,\uv_k \right)
	\leq
		\frac{\eta}{2\gamma^2\eig_r\left(\He\right)}
		\Vnorm{\xv_k-\zv_0}_2^2
\end{align*}
Note that the quantity $\frac{1}{\gamma}\left(\xv_k-\zv_0\right)$ is nothing but the gradient mapping of $p\left(\zv_0,\uv_{k-1}\right)$, hence by \cite[Theorem~1]{Nes13},
\begin{align*}
	\frac{1}{2\gamma}
	\Vnorm{\xv_k-\zv_0}_2^2
	\leq
		p\left(\xv_k,\uv_{k-1}\right)
		-
		p\left(\zv_0,\uv_{k-1}\right)
	\leq
		p\left(\xv_k,\uv_{k-1}\right)
		-
		p\left(\xv_k^\star,\uv_{k-1}\right)	
	\leq 
		\epsilon_{k-1}.
\end{align*}
We thus have
\begin{align*}
	p\left(\zv_0,\uv_k \right)-	p\left(\xv_{k+1}^\star,\uv_k \right)
	\leq
		\ksub\epsilon_{k-1}
	= 
		\frac{\ksub}{1-\rho}\epsilon_k,
\end{align*}
as desired.

\section{Proof of Auxiliary Lemmas}

Recall that at each step of Algorithm~\ref{alg:1}, we wish to solve the following problem:
\begin{align}\label{prob:prox}
	\underset{\xv \in \R^d}{\mbox{minimize}}\left\{q\left(x\right)\triangleq\hF{\xv} + \frac{1}{2\eta} \Mnorm{\xv-\yv + \eta\He^{-1} \vv}{\He}^2\right\}.
\end{align}
The following lemma is a generalization of \cite[Lemma~3]{HKP09} and \cite[Lemma~3]{XZ14} to account for inexactness in the evaluation of the proximal operator, and a Mahalanobis norm.
\begin{lemma}\label{lem:svrg:progress}
Let Assumptions~\ref{asmp:1}--\ref{asmp:2} hold. For any $\yv \in \mathrm{dom}h$, and $\vv \in \R^d$, let $\xv^+$ be an $\epsilon$-optimal solution to  problem~\eqref{prob:prox} in the sense of~\eqref{def:eps:sol},
where $ \eta \in \left(0, \frac{1}{\Lavg} \right]$ is a constant step-size. Denote $\Nvar = \vv - \grad{\fF{}{\yv}}$, then, there exists a vector $\BB{\xi} \in \R^d$ such that $\Mnorm{\BB{\xi}}{\He^{-1}} \leq \sqrt{2\eta\epsilon}$, and it holds for any $\uv \in \R^d$ that
\begin{align*}
	\FF{\xv^+} 
	&\leq
		\FF{\uv}
		-
		\frac{\mu}{2}\Mnorm{\yv-\uv}{\He}^2
		-
		\frac{1}{2\eta} 
		\left( 
			\Mnorm{\xv^+-\uv}{\He}^2
			-
			\Mnorm{\yv-\uv}{\He}^2
		 \right)
	\\
	&\hspace{1.5cm}
		+
		\InP{\Nvar}{\uv-\xv^+}
		+ 
		\InP{\BB{\xi}}{\uv-\xv^+}
		+ 
		\epsilon.
\end{align*}
\begin{proof}
Since $\xv^+$ is an $\epsilon$-optimal solution to problem~\eqref{prob:prox},
by Lemma~\ref{lem:eps:subdiff}, there exists a vector $\BB{\xi}$, satisfying $\Mnorm{\BB{\xi}}{\He^{-1}} \leq \sqrt{2\eta\epsilon}$, such that 
\begin{align}\label{eq:lem:3:proof:1}
	-\frac{1}{\eta}
	\left(	
	\He\left(\xv^+ - \yv\right) 
	+
	\eta\vv
	+
	\BB{\xi}	
	\right)
	\in
	\partial_\epsilon \hF{\xv^+}.
\end{align}
Since $h$ is convex, for any vector $\BB{\zeta} \in	\partial_\epsilon \hF{\xv^+}$, it holds for all $\uv \in \R^{d}$ that 
\begin{align*}
	\hF{\uv} - \hF{\xv^+} \geq \InP{\BB{\zeta}}{\uv-\xv^+}-\epsilon.
\end{align*}
 Therefore, it follows from~\eqref{eq:lem:3:proof:1} that
\begin{align}\label{eq:lem:3:proof:2}
	\InP{
		\vv
		+
		\eta^{-1}
		\He
		\left(\xv^+-\yv\right)
		+
		\eta^{-1}
		\BB{\xi}
	}{\xv^+-\uv}
	+
	\hF{\xv^+}
	-
	\epsilon
	\leq
	\hF{\uv}.
\end{align}
By the smoothness of $f$, we have for $\eta \in \left(0,\frac{1}{\Lavg}\right]$ that
\begin{align*}
	\FF{\xv^+} 
	&\leq
		\fF{}{\yv} + \InP{\grad{\fF{}{\yv}}}{\xv^+-\yv} 
		+
		\frac{1}{2\eta} \Mnorm{\xv^+-\yv}{\He}^2
		+ 
		\hF{\xv^+}
	\\
	&\hspace{0.0cm}
	\leq
		\fF{}{\yv} 
		+		
		\InP{\grad{\fF{}{\yv}}}{\xv^+-\uv}
		+
		\InP{\grad{\fF{}{\yv}}}{\uv-\yv}
		+
		\frac{1}{2\eta} \Mnorm{\xv^+-\yv}{\He}^2
		+ 
		\hF{\xv^+}
	\\
	&\hspace{0.0cm}
	=
		\fF{}{\yv} 
		+		
		\InP{\grad{\fF{}{\yv}}}{\uv-\yv}		
		+
		\InP{\grad{\fF{}{\yv}}}{\xv^+-\uv}
	\\
	&\hspace{1.35cm}	
		+
		\frac{1}{2\eta} \Mnorm{\xv^+-\uv}{\He}^2
		+
		\frac{1}{\eta}
		\InP{\uv-\yv}{\xv^+-\uv}_{\He}		
		+
		\frac{1}{2\eta} \Mnorm{\yv-\uv}{\He}^2
		+ 
		\hF{\xv^+},
\end{align*}
where the last equality follows from adding and subtracting the term $\uv$ and expanding the norm squared. Since $f$ is $\mu$-strongly convex w.r.t the $\He$-norm, we have 
\begin{align*}
	\fF{}{\yv} +\InP{\grad{\fF{}{\yv}}}{\uv-\yv}
	\leq
		\fF{}{\uv}
		-
		\frac{\mu}{2}
		\Mnorm{\yv-\uv}{\He}^2.
\end{align*}
It follows that
\begin{align*}
\FF{\xv^+} 
	&\leq
		\fF{}{\uv}
		-
		\frac{\mu}{2}
		\Mnorm{\yv-\uv}{\He}^2
		+
		\InP{\grad{\fF{}{\yv}}}{\xv^+-\uv}		
	\\
	&\hspace{0.3cm}			
		+
		\frac{1}{2\eta} \Mnorm{\xv^+-\uv}{\He}^2
		+
		\frac{1}{\eta}
		\InP{\uv-\yv}{\xv^+ -\uv}_{\He}		
		+
		\frac{1}{2\eta} \Mnorm{\yv-\uv}{\He}^2
		+ 		
		\hF{\xv^+}
	\\
	&\hspace{0.0cm}
	=
		\fF{}{\uv}
		-
		\frac{\mu}{2}
		\Mnorm{\yv-\uv}{\He}^2
		-
		\frac{1}{2\eta} 
		\left( 
			\Mnorm{\xv^+-\uv}{\He}^2
			-
			\Mnorm{\yv-\uv}{\He}^2
		 \right)
	\\
	&\hspace{.3cm}
		+
		\InP{\Nvar}{\uv-\xv^+}
		+
		\InP{\vv + \eta^{-1}\He\left(\xv^+ -\yv\right) + \eta^{-1}\BB{\xi}}{\xv^+-\uv}
		+ 
		\hF{\xv^+}
		+
		\frac{1}{\eta}
		\InP{\BB{\xi}}{\uv-\xv^+}
	\\
	&\hspace{0.0cm}
	\leq
		\FF{\uv}
			-
		\frac{\mu}{2}
		\Mnorm{\yv-\uv}{\He}^2
		-
		\frac{1}{2\eta} 
		\left( 
		\Mnorm{\xv^+-\uv}{\He}^2
		-
		\Mnorm{\yv-\uv}{\He}^2
		\right)
	\\
	&\hspace{2.5cm}
		+
		\InP{\Nvar}{\uv-\xv^+}
		+
		\frac{1}{\eta}
		\InP{\BB{\xi}}{\uv-\xv^+}
		+\epsilon,
\end{align*}
where the equality follows from adding and subtracting the terms $\vv$ and $\xv^+$ in the first and second inner products, respectively; and the last inequality follows from~\eqref{eq:lem:3:proof:2}. This completes the proof of Lemma~\ref{lem:svrg:progress}.
\end{proof}
\end{lemma}

\begin{lemma}\label{lem:var:grad}
Let Assumptions 1--2 hold. Let $\xv_{k+1}$ be an $\epsilon_k$-optimal solution to the subproblem in Step~7 of Algorithm~\ref{alg:1}. Let $\uv_k$ be any vector in $\R^d$ that is independent of the mini-batch $\mathcal{B}_k$, then, it holds for any $\eta\in \left(0,\frac{1}{\Lavg}\right]$ that 
\begin{align*}
	\EX{\InP{\Nvar_k}{\uv_k-\xv_{k+1}}}
	\leq
		\frac{3\eta}{2}
		\E\Mnorm{\Nvar_k}{\He^{-1}}^2		
		+
		\epsilon_k.
\end{align*}
\begin{proof}
Let $\xv_{k+1}^\star$ be the exact solution to the subproblem in Step~7 of Algorithm~\ref{alg:1}, i.e, $\epsilon_k=0$. Since the objective function $q$ defined in~\eqref{prob:prox} is $\frac{1}{\eta}$-strongly convex w.r.t the $\He$-norm, we have
\begin{align*}
	\frac{1}{2\eta}\Mnorm{\xv_{k+1}-\xv_{k+1}^\star}{\He}^2 \leq q\left(\xv_{k+1}\right)-q\left(\xv_{k+1}^\star\right) \leq \epsilon_k.
\end{align*}
Therefore, one can write $\xv_{k+1}$ as $\xv_{k+1}=\xv_{k+1}^\star+\BB{\zeta}_k$ for some vector $\BB{\zeta}_k \in \R^d$ satisfying $\Mnorm{\BB{\zeta}_k}{\He}\leq\sqrt{2\eta\epsilon_k}$. 
If we define the following virtual iterate:
\begin{align*}
	\bar{\yv}_k = \Prox{\eta h}{\He}{\yv_k - \eta\He^{-1}\grad{\fF{}{\yv_k}}},
\end{align*}
then
\begin{align}\label{eq:lem:var:grad:proof:1}
	\InP{\Nvar_k}{\uv_k-\xv_{k+1}}
	&= 
		\InP{\Nvar_k }{\uv_k-\bar{\yv}_k}
		+
		\InP{\Nvar_k }{\bar{\yv}_k-\xv_{k+1}}
	\nonumber\\
	&\leq
		\InP{\Nvar_k }{\uv_k-\bar{\yv}_k}
		+
		\Mnorm{\Nvar_k}{\He^{-1}}
		\Mnorm{\bar{\yv}_k-\xv_{k+1}}{\He},
\end{align}
where the last step follows from Cauchy's inequality.
Note that 
\begin{align*}
	\Mnorm{\bar{\yv}_k-\xv_{k+1}}{\He}
	&=
	\Mnorm{
		\Prox{\eta h}{\He}{\yv_k - \eta \He^{-1}\grad{\fF{}{\yv_k}}}
		-
		\Prox{\eta h}{\He}{\yv_k - \eta \He^{-1}\vv_k} 
		-
		\BB{\zeta_k}
	}{\He}
	\\
	&\leq
	\Mnorm{
			\Prox{\eta h}{\He}{\yv_k - \eta \He^{-1}\grad{\fF{}{\yv_k}}}
			-
			\Prox{\eta h}{\He}{\yv_k - \eta \He^{-1}\vv_k} 
	}{\He}
	+
	\Mnorm{\BB{\zeta_k}}{\He}
	\nonumber\\
	&\leq
		\eta\Mnorm{\Nvar_k}{\He^{-1}}
		+
		\sqrt{2\eta\epsilon_k},
\end{align*}
where the first inequality follows from the triangle inequality and the last one follows from Property~\ref{lem:prox}. 
Therefore, it holds that
\begin{align*}
	\Mnorm{\Nvar_k}{\He^{-1}}
	\Mnorm{\bar{\yv}_k-\xv_{k+1}}{\He}
	&\leq 
		\eta
		\Mnorm{\Nvar_k}{\He^{-1}}^2
		+
		\sqrt{2\eta\epsilon_k}		
		\Mnorm{\Nvar_k}{\He^{-1}}
	\nonumber\\
	&\leq 
		\frac{3\eta}{2}
		\Mnorm{\Nvar_k}{\He^{-1}}^2		
		+
	\epsilon_k.
\end{align*}
Denote the filtration $\mathcal{F}_k$ by $\mathcal{F}_k=\{\BB{f}_0, \BB{f}_1,\ldots, \BB{f}_{k}\}$, where $\BB{f}_{i}$ denotes all the randomness incurring at time $i$ for all $i\in\{1,\ldots,k\}$. Note that the triple $\left(\xv_k,\yv_k,\zv_k\right)$ depends on $\mathcal{F}_{k-1}$, but not on $\BB{f}_{k}$, we thus have 
\begin{align*}
	\E_{\mathcal{F}_k}{\InP{\Nvar_k }{\bar{\yv}_k-\uv_k}}
	=
		\E_{\mathcal{F}_{k-1}}{
		\E_{\BB{f}_{k}}{\InP{\Nvar_k }{\bar{\yv}_k-\uv_k|\mathcal{F}_{k-1}}}
		}
	=
		\E_{\mathcal{F}_{k-1}}{
		\InP{\E_{\BB{f}_{k}}{\Nvar_k}}{\bar{\yv}_k-\uv_k|\mathcal{F}_{k-1}}
	}
	=0.
\end{align*}
Taking the expectation on both sides of~\eqref{eq:lem:var:grad:proof:1} completes the proof.
\end{proof}
\end{lemma}

\begin{lemma}[Bounding Variance]\label{lem:bouding:var}
Assume that the indices in the mini-batch $\mathcal{B}_k$ are sampled independently from $\{1,\ldots,n\}$ with probabilities $P_i=\frac{L_i}{n\Lavg}$, then conditioned on $\xv_k$ and $\yv_k$, it holds that
\begin{align*}
	\E{\Mnorm{\Nvar_k}{\He^{-1}}^2}
	\leq	
		\frac{2\Lavg^2}{b}	
		\Mnorm{\xv_k-\yv_k}{\He}^2
		+
		\frac{8\Lavg}{b}
		\left(	
		\FF{\xv_k} - \FF{\xv^\star}
		+
		\FF{\xTd_s} - \FF{\xv^\star}\right).
\end{align*}
\begin{proof}
We follow the original proofs in~\cite{XZ14,Nit14} with a few modifications to take the Mahalanobis norm into account. For any $i\in\{1,\ldots,n\}$, consider the function $\phi_i$ defined by
\begin{align*}
	\phi_i\left(\xv\right)=\fF{i}{\xv}-\fF{i}{\xv^\star}-\InP{\grad{\fF{i}{\xv^\star}}}{\xv-\xv^\star}.
\end{align*}
Then, $\phi_i\left(\xv^\star\right)=\min_{\xv} \phi_i\left(\xv\right)$ since $\grad{\phi_i\left(\xv^\star\right)}=0$ and $\phi_i$ is convex. 
It can be checked that $\grad{\phi_i}$ is $L_i$-Lipschitz, we thus have 
\begin{align*}
	0=\phi_i\left(\xv^\star\right)
	&\leq
		\min_{\eta}
		\phi_i\left(\xv-\eta\He^{-1}\grad{\phi_i\left(\xv\right)}\right)
	\nonumber\\
	&\leq
		\min_{\eta}
		\phi_i\left(\xv\right)
		-
		\eta
		\Mnorm{\grad{\phi_i\left(\xv\right)}}{\He^{-1}}^2
		+
		\frac{\eta^2L_i}{2}
		\Mnorm{\grad{\phi_i\left(\xv\right)}}{\He^{-1}}^2
	\nonumber\\
	&=
		\phi_i\left(\xv\right)
		-
		\frac{1}{2L_i}
		\Mnorm{\grad{\phi_i\left(\xv\right)}}{\He^{-1}}^2,
\end{align*}
which is equivalent to
\begin{align*}
	\Mnorm{\grad{\fF{i}{\xv}}-\grad{\fF{i}{\xv^\star}}}{\He^{-1}}^2
	\leq
		2L_i\left(\fF{i}{\xv}-\fF{i}{\xv^\star}-\InP{\grad{\fF{i}{\xv^\star}}}{\xv-\xv^\star}\right).
\end{align*}
Multiplying both sides of the previous inequality by $\frac{\Lavg}{L_i}$ and averaging from $1,\ldots,n$ gives
\begin{align*}
	\frac{1}{n}
	\sum_{i=1}^n
	\frac{\Lavg}{L_i}
	\Mnorm{	\grad{\fF{i}{\xv}}-\grad{\fF{i}{\xv^\star}}}{\He^{-1}}^2
	\leq
		2\Lavg\left(\fF{}{\xv}-\fF{}{\xv^\star}\right)
		-
		\InP{\grad{\fF{}{\xv^\star}}}{\xv-\xv^\star}.
\end{align*}
Since $\FF{\xv}=\fF{}{\xv}+\hF{\xv}$,  by the optimality of $\xv^\star$, $-\grad{\fF{}{\xv^\star}}\in\partial \hF{\xv^\star}$. Therefore, it follows from the convexity of $h$ that
$
	-\InP{\grad{\fF{}{\xv^\star}}}{\xv-\xv^\star}
	\leq
		\hF{\xv}-\hF{\xv^\star},
$
which implies that
\begin{align}\label{eq:sum:grad}
	\frac{1}{n}
	\sum_{i=1}^n
	\frac{\Lavg}{L_i}
	\Mnorm{\grad{\fF{i}{\xv}}-\grad{\fF{i}{\xv^\star}}}{\He^{-1}}^2
	\leq
	2\Lavg \left(\FF{\xv}-\FF{\xv^\star}\right).
\end{align}

We are now ready to bound $\E{\Mnorm{\Nvar_k}{\He^{-1}}^2}$, we have
\begin{align*}
	\E\Mnorm{\vv_k-\grad{\fF{}{\yv_k}}}{\He^{-1}}^2
	&=
		\frac{1}{b^2}
		\E\bigg\{ 
			\sum_{i\in\mathcal{B}_k}
			\Mnorm{
				\frac{
						\grad{\fF{i}{\yv_k}} 								
						-
						\grad{\fF{i}{\xTd_{s}}} 				
					}{L_i/\Lavg}
				-
				\left(
					\grad{\fF{}{\yv_k}}
					-
					\grad{\fF{}{\xTd_{s}}}
				\right)
			}{\He^{-1}}^2		
		\bigg\}
	\nonumber\\
	&\leq
		\frac{\Lavg}{b}
		\E
		\frac{1}{L_i}
			\Mnorm{
				\grad{\fF{i}{\yv_k}} 
				- 
				\grad{\fF{i}{\xTd_{s}}} 			
		}{\He^{-1}}^2
	\nonumber\\
	&\leq 
		\frac{2\Lavg}{b}
		\E
		\frac{1}{L_i}
		\Mnorm{
			\grad{\fF{i}{\yv_k}} 
			- 
			\grad{\fF{i}{\xv_k}} 			
		}{\He^{-1}}^2		
	\nonumber\\
	&\hspace{0.5cm}
		+
		\frac{4\Lavg}{b}
		\E
		\frac{1}{L_i}
		\Mnorm{
			\grad{\fF{i}{\xv_k}} 
			- 
			\grad{\fF{i}{\xv^\star}} 			
		}{\He^{-1}}^2
	\nonumber\\
	&\hspace{0.5cm}
		+ 
		\frac{4\Lavg}{b}	
		\E
			\frac{1}{L_i}
			\Mnorm{
			\grad{\fF{i}{\xTd_{s}}} 			
			-
			\grad{\fF{i}{\xv^\star}} 			
		}{\He^{-1}}^2		
	\nonumber\\
	&\leq			
		\frac{2\Lavg^2}{b}
		\Mnorm{\xv_k-\yv_k}{\He}^2		
		+
		\frac{8\Lavg}{b}	
		\left(
			\FF{\xv_k} - \FF{\xv^\star}
			+
			\FF{\xTd_s} - \FF{\xv^\star}
		\right),		
\end{align*}
where the first inequality follows from the fact that $\E\Vnorm{X-\E X}^2\leq \E\Vnorm{X}^2$; we have used the inequality $\Mnorm{\xv+\yv}{\He^{-1}}^2\leq 2 \Mnorm{\xv}{\He^{-1}}^2 + 2\Mnorm{\xv}{\He^{-1}}^2$ to derive the second inequality; and the last step follows by evaluating the expectations, and using Assumption~\ref{asmp:1} and~\eqref{eq:sum:grad}. This completes the proof of Lemma~\ref{lem:bouding:var}.
\end{proof}
\end{lemma}

The following definition of  $\epsilon$-subgradients is very useful for analyzing \replaced{how ineaxact proximal evaluations affect the convergence of the algorith.}{the effect of inexactness in the proximal evaluations on the convergence of the algorithm.}  
\begin{definition}[$\epsilon$-Subgradients \cite{BNO03}]
Given a convex function $f: \R^d\to \R$  and a positive constant $\epsilon$, we say that a vector $\BB{\zeta}\in\R^d$  is an $\epsilon$-subgradient of $f$ at $\xv\in\R^d$ if
\begin{align*}
	f\left(\zv\right)\geq f\left(\xv\right) + \InP{\BB{\zeta}}{\zv-\xv}-\epsilon, \quad \forall \zv\in\R^d.
\end{align*}
The set of all $\epsilon$-subgradients of $f$ at $\xv$ is called the $\epsilon$-subdifferential of $f$ at $\xv$, and is denoted by $\partial_\epsilon f\left(\xv\right)$.
\end{definition}

The following lemma characterizes the property of the $\epsilon$-differential of the function $h$, where its proof for the case of the Euclidean norm  can be found in \cite[Lemma~2]{SRB11}. We provide the proof here for completeness.
\begin{lemma}\label{lem:eps:subdiff}
If $\xv^+$ is an $\epsilon$-optimal solution to problem~\eqref{prob:prox} in the sense of \eqref{def:eps:sol}, then there exists a vector $\BB{\xi}\in\R^d$ such that $\Mnorm{\BB{\xi}}{\He^{-1}} \leq \sqrt{2\eta\epsilon}$ and 
\begin{align*}	
	-\frac{1}{\eta}
	\left(	
		\He\left(\xv^+ - \yv\right) 
		+
		\eta\vv
		+
		\BB{\xi}	
	\right)
	\in
	 \partial_\epsilon \hF{\xv^+}.
\end{align*}
\begin{proof}
We start by noting that $\B{0}\in\partial_\epsilon q\left(\xv^+\right)$ when $\xv^+$ is an $\epsilon$-optimal solution to problem~\eqref{prob:prox}. Consider the function $f\left(\xv^+\right)=\frac{1}{2\eta}\Mnorm{\xv^+-\tv}{\He}^2$, then it can be verified that
\begin{align*}
	\partial_\epsilon f\left(\xv^+\right)
	&=
		\big\{\BB{\zeta}\in\R^d\,
		\big| \,
			\frac{1}{2\eta}
			\Mnorm{\He\left(\xv^+-\tv \right)-\eta\BB{\zeta}}{\He^{-1}}^2
		\leq
		\epsilon
		\}
	\\
	&=
		\big\{\BB{\zeta}\in\R^d, 
		\BB{\zeta}=\frac{1}{\eta}\He\left(\xv^+-\tv \right)+\frac{1}{\eta}\BB{\xi}
		\,\big| \,
			\frac{1}{2\eta}
			\Mnorm{\BB{\xi}}{\He^{-1}}^2
		\leq
		\epsilon
		\}.			
\end{align*}
We have for convex functions $f_1,f_2$ that $\partial_\epsilon \left(f_1+f_2\right)\left(\xv\right)\subset \partial_\epsilon f_1\left(\xv\right)+\partial_\epsilon f_2\left(\xv\right)$ \cite{BNO03}. Therefore, if we let $f_1=f$, $f_2=h$, and $\tv=\yv - \eta\He^{-1} \vv$, then $\B{0}\in\partial_\epsilon q\left(\xv^+\right)$. Since $\partial_\epsilon q\left(\xv^+\right)\subset \partial_\epsilon f\left(\xv^+\right)+\partial_\epsilon h\left(\xv^+\right)$, it follows that $\B{0}$ must be a sum of an element of $\partial_\epsilon f\left(\xv^+\right)$ and an element of $\partial_\epsilon h\left(\xv^+\right)$. Thus, there is a vector $\BB{\xi}\in \R^d$ such that
\begin{align*}
	-\frac{1}{\eta}\He\left(\xv^+-\tv \right)-\frac{1}{\eta}\BB{\xi}
	\in 
	\partial_\epsilon h\left(\xv^+\right)
	\quad \text{with}\quad \Mnorm{\BB{\xi}}{\He^{-1}} \leq \sqrt{2\eta\epsilon},
\end{align*}
completing the proof.
\end{proof}
\end{lemma}

\begin{lemma}[{\cite[Lemma~1]{SRB11}}]\label{lem:iter:recurs}
Assume that the nonnegative sequence ${u_k}$ satisfies the following recursion for all $k\geq 1$:
\begin{align*}
	u_k^2 \leq S_k + \sum_{i=1}^{k}\alpha_t u_t,
\end{align*}
where $\{S_k\}$ is an increasing sequence, $S_0\geq u_0^2$, and $\alpha_t\geq 0$ for all $t$. Then, for all $k \geq 1$, then
\begin{align*}
	u_k 
	\leq
		\frac{1}{2}
		\sum_{t=1}^{k}\alpha_t
		+
		\left(
			S_k
			+
			\left(\frac{1}{2}\sum_{t=1}^{k}\alpha_t\right)^2
		\right)^{\frac{1}{2}}. 
\end{align*}
\end{lemma}

\end{document}